\documentclass[journal]{IEEEtran}
\usepackage[noadjust]{cite}
\pdfoutput=1

\usepackage{graphicx}

\usepackage{amsmath}
\usepackage{array}

\ifCLASSOPTIONcompsoc
  \usepackage[caption=false,font=normalsize,labelfont=sf,textfont=sf]{subfig}
\else
  \usepackage[caption=false,font=footnotesize]{subfig}
\fi
\usepackage{url}
\usepackage{amssymb}

\usepackage{algorithm}
\usepackage{algorithmicx}
\usepackage[noend]{algpseudocode}
\newcommand*\Let[2]{\State #1 $\gets$ #2}
\algrenewcommand\algorithmicrequire{\textbf{Input:}}
\algrenewcommand\algorithmicensure{\textbf{Output:}}

\usepackage[normalem]{ulem}
\usepackage{tikz}
\usepackage{color}
\definecolor{mygreen}{rgb}{0 0.7 0}
\definecolor{mygray}{rgb}{0.6 0.6 0.6}
\definecolor{myblue}{rgb}{0.7 0.7 1}
\definecolor{mypink}{rgb}{1 0 1}

\usepackage{multirow}
\usepackage[pdftex]{mathtools}
\def\block(#1,#2)#3{\multicolumn{#2}{c}{\multirow{#1}{*}{$ #3 $}}}

\begin{document}

\title{Sensor Selection With Cost Constraints for Dynamically Relevant Bases}

\author{Emily~Clark\thanks{E. Clark is with the Department of Physics at
the University of Washington, Seattle, WA, 98195-1560 USA email:
eclark7@uw.edu},
        J.~Nathan~Kutz\thanks{J. N. Kutz is with the Department of Applied
Mathematics at the University of Washington, Seattle, WA, 98195-3925 
USA email: kutz@uw.edu},~\IEEEmembership{Member,~IEEE},
        Steven~L.~Brunton\thanks{S. Brunton is with the Department of Mechanical Engineering
at the University of Washington, Seattle, WA, 98195-2600 USA email:
sbrunton@uw.edu},~\IEEEmembership{Senior Member,~IEEE}}

\markboth{arXiv Submission}%
{Clark \MakeLowercase{\textit{et al.}}: Sensor Placement with Cost Constraints}

\maketitle

\begin{abstract}
We consider cost-constrained sparse sensor selection for full-state reconstruction, applying a well-known greedy algorithm to dynamical systems for which the usual singular value decomposition (SVD) basis may not be available or preferred. 
We apply the cost-modified, column-pivoted QR decomposition to a physically relevant basis\textemdash the pivots correspond to sensor locations, and these locations are penalized with a heterogeneous cost function. 
In considering different bases, we are able to account for the dynamics of the particular system, yielding sensor arrays that are nearly Pareto optimal in sensor cost and performance in the chosen basis. 
This flexibility extends our framework to include actuation and dynamic estimation, and to select sensors without training data.  
We provide three examples from the physical and engineering sciences and evaluate sensor selection in three dynamically relevant bases: truncated balanced modes for control systems, dynamic mode decomposition (DMD) modes, and a basis of analytic modes. 
We find that these bases all yield effective sensor arrays and reconstructions for their respective systems. 
When possible, we compare to results using an SVD basis and evaluate tradeoffs between methods.
\end{abstract}



%



\section{Introduction}
\label{sec:intro}
\IEEEPARstart{M}{odal} decompositions and model reduction are fundamental to many fields of science and engineering~\cite{Taira2017aiaa,Taira2019aiaa,Brunton2019book}. 
Historically, such decompositions were based on 
special functions, such as Fourier modes, Legendre polynomials, and Bessel functions. 
The development in the late 1960s~\cite{golub1971singular} of an efficient numerical procedure for computing the {\em singular value decomposition} (SVD) allowed for data-driven modal decompositions that far outperformed generalized functional techniques. 
The SVD is now the dominant paradigm for model reduction in many applications, and for good reason\textemdash for a given number of modes, it is guaranteed to capture the greatest amount of the system's energy (variance)~\cite{eckart1936approximation}. 
However, ordering modes based on energy does not result in the most effective basis for every task.  
For control systems, modes obtained via balanced truncation~\cite{moore1981principal} are often more effective than SVD modes. 
Similarly, the dynamic mode decomposition (DMD)~\cite{schmid2010dynamic,rowley2009spectral,tu2014dynamic,Kutz2016book} provides an alternative basis that simultaneously captures spatial correlation and evolution in time.  
%
%
Sparse sensor optimization is often based on SVD modes~\cite{manohar2018data,clark2018greedy}, although this approach has been previously extended to include these alternative bases~\cite{manohar2019optimized,manohar2018optimal}.   
The goal of this work is to develop a greedy sensor selection algorithm that incorporates a heterogeneous cost on sensor location and apply it to these alternative bases.  
Accounting for a heterogeneous cost on sensor location is a critical practical aspect of sensor selection, and our unified framework results in sensor arrays that are approximately Pareto optimal in cost and performance.  

The full state-space of a high-dimensional, spatio-temporal system may be prohibitively expensive to measure, thus requiring significant measurement downsampling.   Sparse sensor selection provides a mathematical framework for how to best downsample in many applications, including reduced order modeling~\cite{willcox2006unsteady}, scientific experiments~\cite{antchev2010totem}, control theory~\cite{lim1992method}, climate science~\cite{yang2010eof}, and industry~\cite{liu2015entropy,Manohar2018jms}. We formulate sensor selection as an interpolation problem, allowing us to adopt methods from nonlinear model order reduction, such as the empirical and discrete empirical interpolation methods (DEIM)~\cite{barrault2004empirical, chaturantabut2010nonlinear} and Q-DEIM~\cite{drmac2016new}. Manohar et al~\cite{manohar2018data} thoroughly demonstrate the use of column-pivoted QR for sensor selection for reconstruction in an SVD basis, and in~\cite{clark2018greedy} the authors expand upon this, applying it to a randomized reduced basis~\cite{liberty2007randomized, halko2011finding} and modifying the QR algorithm to account for a heterogeneous cost on sensor location.

There are other sensor selection paradigms, most notably methods that exploit the submodular structure of certain performance metrics to develop greedy algorithms with performance guarantees, e.g.~\cite{krause2008near, ranieri2014near}. Some of these algorithms can also be adapted to account for a non-heterogeneous sensor cost~\cite{krause2006near, krause2007near, krause2008robust, krause2014submodular}. However, the greedy approximate determinant maximization provided by column-pivoted QR yields good stability and performance guarantees, and is easy to implement. Therefore, we choose the QR-based sensor selection method, even though the determinant is not necessarily submodular.

The main contribution of this work is applying the cost-constrained column-pivoted QR algorithm to new, dynamically relevant bases. This extension of the method increases the practicality of sensor selection by both accounting for sensor cost and using a basis that is matched to the task at hand (i.e., reconstruction, estimation, or control). 
First, we demonstrate balanced sensor and actuator selection for control systems, following the method of Manohar et al.~\cite{manohar2018optimal}, which we modify to include a cost function. Next, we perform sensor selection on a dynamic mode decomposition (DMD) basis~\cite{schmid2010dynamic}. 
The DMD gives approximate equations of motion for the system, allowing us to construct a Kalman filter that outperforms interpolative reconstructions in the presence of noise. Finally, we apply our method to universal basis modes, such as Fourier modes, Bessel functions, or a polynomial basis. 
This is useful for systems that admit a universal basis and for which the equations of motion are known and solvable for a modal decomposition. 
 In this way, it is possible to perform principled sensor selection without training data. In all three applications, we show that the modified column-pivoted QR decomposition yields sensor (and actuator) arrays with accurate reconstructions, while accounting for costs on sensor (and actuator) location. For the DMD and analytic mode examples, we also provide comparisons with sensors selected using an SVD basis. We find that the SVD basis and sensors often lead to better reconstructions, but the dynamically relevant bases have other advantages\textemdash information about the dynamics in the case of DMD modes, and no need for high-fidelity training data with analytically calculated modes.

The remainder of the paper is organized as follows: In Section \ref{sec:algorithm}, we present QR pivoting for sensor selection and the modification to incorporate a cost function. 
In Section \ref{sec:control}, we describe sensor and actuator selection for control systems, including an example from physics. Section \ref{sec:DMD} considers cost-constrained sensor selection for DMD modes, with an example from climate science.  Section \ref{sec:analytic} explores sensor selection with cost for universal basis modes, demonstrated on another example from physics. We present conclusions in Section \ref{sec:conclusion}.

\section{Principled sparse sensor selection}
\label{sec:algorithm}
\subsection{Problem Formulation}
\label{subsec:problem}
We consider a high-dimensional system, such as the discretization of a PDE or a set of high-definition photographs. We wish to sparsely sample the state of the system and then reconstruct the full state from these sparse measurements. The goal is to determine the optimal sampling points that yield the best reconstructions for a given number of sensors.  
We will train these sensor configurations based on a set of $m$ full-state snapshots ${\bf x}_i\in\mathbb{R}^n$, which we gather into a data matrix
\begin{align}
{\bf X} = \left[\begin{array}{cccc}
{\bf x}_1 & {\bf x}_2 & \cdots & {\bf x}_m
\end{array}\right].
\end{align}

We select a basis ${\bf\Psi}\in\mathbb{R}^{n\times r}$ in which to represent the snapshots, such that
\begin{align}\label{Eq:Representation}
{\bf x} &= {\bf\Psi}{\bf a},
\end{align}
where ${\bf a}\in\mathbb{R}^r$. Though it is not a requirement, usually $r\ll n$, so that ${\bf\Psi}$ is a \emph{reduced} basis.    

The goal of sensor selection is to learn a selection operator ${\bf C}$ that extracts $p\ll n$ point measurements:
\begin{align}
{\bf Y} = {\bf C}{\bf X},
\label{eq:measure}
\end{align}
where ${\bf Y} \in\mathbb{R}^{p\times m}$ are the measurements. 
The rows of ${\bf C}\in\mathbb{R}^{p\times n}$ are unit vectors, picking out rows of ${\bf X}$ (which correspond to locations in space) to measure.

In the ${\bf\Psi}$ basis, Eq. \ref{eq:measure} becomes
\begin{align}
{\bf Y} = {\bf C}{\bf\Psi}{\bf a}
= {\bf\Theta}{\bf a}.
\label{eq:sense}
\end{align}
We obtain the full-state reconstructions by taking the minimum-norm, least squares solution,
\begin{align}
\label{eq:recona}
{\bf \hat{a}} &= {\bf\Theta}^\dagger{\bf Y}\\
\label{eq:reconx}
{\bf\hat{X}} &= {\bf\Psi}{\bf\hat{a}},
\end{align}
with ${\bf\Theta}^\dagger$ being the pseudoinverse of ${\bf\Theta}$.

The fractional reconstruction error is: 
\begin{align}
E = \frac{||{\bf X}-{\bf\hat{X}}||_F}{||{\bf X}||_F},
\label{eq:error}
\end{align}
where $||\cdot ||_F$ denotes the Frobenius norm.

The sensor selection problem is to determine the optimal set of indices $\{J\} \in [1, n]$, $|J| = p$, at which to measure the system, so as to minimize the reconstruction error. In other words, we build ${\bf C}$ by setting its rows to ${\bf e}_{\{J\}}^T$, where ${\bf e}_j$ is the unit vector in the $j^{\text{th}}$ direction. This is a combinatorially hard problem, but as described in the introduction and below in subsection \ref{subsec:QR}, there are greedy and heuristic methods that yield near-optimal solutions.

We further extend the problem to include a heterogeneous cost function $\boldsymbol{\eta}$ on the set of candidate sensor locations. We then wish to determine the set of indices $\{J\}$ that is Pareto optimal in reconstruction error and total sensor cost $c$, where
\begin{align}
c = \sum_{j\in\{J\}} \eta_j.
\end{align}

Since ${\bf\Theta} = {\bf C}{\bf\Psi}$, the sensor selection problem is equivalent to choosing the most informative rows of ${\bf\Psi}$, and as such, the performance of the sensor selection algorithms and the reconstructions depends greatly on the choice of ${\bf\Psi}$. The basis could be derived from full-state training data, e.g., ${\bf\Psi}$ is usually taken to be the first $r$ left singular vectors of ${\bf X}$~\cite{manohar2018data}. This basis choice often gives excellent results, but it is not always optimal for the task at hand. 
For a control system, the SVD provides no information about controllability and observability.  
Similarly, the SVD does not provide a model for how modes evolve in time, and thus does not admit Kalman filtering as an alternative reconstruction method. 
Finally, the SVD requires full-state training data, which may not be available. Accordingly, the main contribution of this paper will be to extend a previously developed cost-constrained sensor selection algorithm~\cite{clark2018greedy} to new bases tailored for several kinds of dynamical systems.

\subsection{Column-pivoted QR}
\label{subsec:QR}
Following~\cite{drmac2016new}, we use the QR decomposition with column pivoting for sensor selection. The QR decomposition is a greedy method for approximate determinant maximization. The pivoted QR decomposition was introduced by Businger and Golub~\cite{businger1965linear} as a method of solving least squares problems, and it was first employed for sensor selection by Drma{\v{c}} and Gugercin~\cite{drmac2016new}, who showed that it provides improved error bounds over the discrete empirical interpolation method.

The pivoted QR decomposition for some matrix ${\bf V}\in\mathbb{R}^{r\times n}$, $r\le n$ is given by:
\begin{equation}
\label{eq:QR}
{\bf V}{\bf P} = {\bf Q}{\bf R},
\end{equation}
where ${\bf Q}\in\mathbb{R}^{r\times r}$ is unitary, ${\bf R}\in\mathbb{R}^{r\times n}$ is upper triangular, and ${\bf P}\in\mathbb{R}^{n\times n}$ is a column permutation matrix.

To perform the decomposition, we first initialize the permutation matrix to ${\bf P}_0 = \mathbb{I}_n$. Then we select the column ${\bf v}_j$ of ${\bf V}$ with the largest norm and swap it with the first column of ${\bf V}$ via the permutation matrix:
\begin{equation}
{\bf P}_1 = \left(\begin{array}{cccccccc}
{\bf e}_j & {\bf e}_2 & \cdots & {\bf e}_{j-1} & {\bf e}_1 & {\bf e}_{j+1} & \cdots & {\bf e}_n
\end{array}\right).
\label{eq:permutation}
\end{equation}
Next we apply a Householder transformation ${\bf Q}_1$ such that
\begin{equation}
{\bf Q}_1{\bf V}{\bf P}_1 = \left(\begin{array}{cccc}
||{\bf v}_j|| & * & * & \cdots\\[4pt]
0\\[4pt]
0 & & {\bf V}'\\
\vdots
\end{array}\right).
\label{eq:Householder}
\end{equation}
We repeat the process on the submatrix ${\bf V}'$, i.e. identify the column $j'$ with the largest norm, swap it with the first column of ${\bf V}'$, and transform it so that its subdiagonal elements are zero:
\begin{equation}
\label{eq:QR_step2}
{\bf Q}_2{\bf Q}_1{\bf V}{\bf P}_2 = \left(\begin{array}{cccc}
||{\bf v}_j|| & \star &\star & \cdots\\[4pt]
0 & ||{\bf v}'_{j'}|| & \star & \cdots\\[4pt]
0 & 0 & \block(2,2){{\bf V}''}\\
\vdots & \vdots
\end{array}\right).
\end{equation}
We continue to iterate until the right-hand side is upper triangular and the decomposition is achieved:
\begin{equation}
{\bf Q}_r \cdots {\bf Q}_1{\bf V}{\bf P}_r = {\bf R}.
\end{equation}
Thus ${\bf P}$ from Eq. \ref{eq:QR} is ${\bf P}_r$, and ${\bf Q} = {\bf Q}_1^*\cdots{\bf Q}_r^*$.

Note that the $k^{\text{th}}$ Householder transformation is acting on a vector of length $r-k+1$, and so ${\bf Q}_k$ is augmented by the appropriate identity matrix. Also see Algorithm \ref{alg:qr}.

\begin{algorithm}[b]
  \caption[Caption]{\label{alg:qr} QR pivoting.}
  \begin{algorithmic}[1]
    \Require{matrix ${\bf V}\in\mathbb{R}^{r\times n}$, number of sensors $p\le r$}
    \Ensure{Partial permutation matrix ${\bf P}$}
    \Let{${\bf P}$}{$\mathbb{I}_n$}
    \For{$k = 1,\ldots,p$}
    \Let{$j_k$}{$\arg\max\limits_{j\ge k}||{\bf V}_{k:r,j}||$}
    \State \texttt{Swap}$({\bf P}_{:,k},{\bf P}_{:,j_k})$.
    \State Calculate Householder reflection ${\bf Q}_k$ such that ${\bf Q}_k \left(\begin{array}{cccc}
	V_{k,j_k} & V_{k+1,j_k} & \cdots & V_{r,j_k}\end{array}\right)^T = \left(\begin{array}{cccc}
	||{\bf V}_{k:r,j_k}|| & 0 & \cdots & 0 \end{array}\right)^T$.
    \Let{${\bf V}$}{$\text{diag}(\mathbb{I}_{k-1},{\bf Q}_k){\bf V}{\bf P}$}
    \EndFor
    \Return{${\bf P}_{:,1:p}$}
  \end{algorithmic}
\end{algorithm}

This algorithm leads to a diagonal dominance structure in ${\bf R}$:
\begin{equation}
|r_{ii}|^2 \ge \sum_{j=1}^k |r_{jk}|^2, \hspace{24pt} 1\le i\le k\le r.
\end{equation}
Because the Householder transformations are unitary, at every step this method greedily maximizes the determinant of ${\bf R}_{1:k,1:k}$, the $k\times k$ submatrix that has been made upper triangular at the $k^{\text{th}}$ step:
\begin{equation}
|\det {\bf R}_{1:k,1:k}| = \prod_{i=1}^k |r_{ii}| = \prod_{i=1}^k \sigma_i,
\end{equation}
or equivalently, it maximizes the volume of the submatrix.

For sensor selection, we perform column-pivoted QR on ${\bf\Psi}^*$, and for $p$ sensors, take ${\bf C}^T = {\bf P}_{:,1:p}$, where ${\bf P}_{:,1:p}$ are the first $p$ columns of ${\bf P}$. Thus we select sensors to approximately maximize the determinant of the measurement matrix ${\bf\Theta}$, a method also known as D-optimal experiment design.

Notice that our problem setup and Algorithm \ref{alg:qr} do not require ${\bf\Psi}$ to be orthogonal, which is what allows us to apply QR pivoting to any appropriate basis for sensor selection. 

\subsection{Cost-constrained sensor selection}
\label{subsec:CCQR}
\begin{algorithm}[b]
  \caption[Caption]{\label{alg:ccqr} QR pivoting with cost function.}
  \begin{algorithmic}[1]
    \Require{matrix ${\bf V}\in\mathbb{R}^{r\times n}$, cost function $\boldsymbol{\eta}\in\mathbb{R}^n$, weighting $\gamma$, number of sensors $p\le r$}
    \Ensure{Partial permutation matrix ${\bf P}$}
    \Let{${\bf P}$}{$\mathbb{I}_n$}
    \For{$k = 1,\ldots,p$}
    \Let{$j_k$}{$\arg\max\limits_{j\ge k}||{\bf V}_{k:r,j}|| - \gamma\eta_j$}
    \State \texttt{Swap}$({\bf P}_{:,k},{\bf P}_{:,j_k})$.
    \State Calculate Householder reflection ${\bf Q}_k$ such that ${\bf Q}_k \left(\begin{array}{cccc}
	V_{k,j_k} & V_{k+1,j_k} & \cdots & V_{r,j_k}\end{array}\right)^T = \left(\begin{array}{cccc}
	||{\bf V}_{k:r,j_k}|| & 0 & \cdots & 0 \end{array}\right)^T$.
    \Let{${\bf V}$}{$\text{diag}(\mathbb{I}_{k-1},{\bf Q}_k){\bf V}{\bf P}$}
    \Let{$\boldsymbol{\eta}$}{${\bf P}\boldsymbol{\eta}$}
    \EndFor
    \Return{${\bf P}_{:,1:p}$}
  \end{algorithmic}
\end{algorithm}

We make a straightforward modification to Algorithm \ref{alg:qr} to account for an additional cost on sensor location. Doing so allows us to select informative sensors at cheaper locations, albeit at some loss of reconstruction accuracy. It is possible to tune the weighting on the cost function to achieve the desired balance between sensor cost and performance.

Assume there is some non-negative cost function $\boldsymbol{\eta}\in\mathbb{R}^n$ and select a scalar weighting $\gamma$. Then perform a modified version of the column-pivoted QR decomposition on ${\bf\Psi}^*$, where at every iteration, select as a pivot the column that maximizes the following quantity~\cite{clark2018greedy}:
\begin{equation}
||\left({\bf\Psi}^*\right)^{(k)}_j|| - \gamma\eta_j.
\end{equation}
Also see Algorithm \ref{alg:ccqr}.

Although the determinant maximization guarantees of column-pivoted QR no longer apply to this modified version, cost-constrained QR balances determinant maximization with cost savings. The weighting $\gamma$ can be varied to trace out a cost versus error curve, with $\gamma=0$ being identical to Algorithm \ref{alg:qr} and a very large value of $\gamma$ relative to the given system and cost function approaching the minimum-cost sensor array. We will show in the following sections that this method often approximately traces out the Pareto curve of simultaneously minimum reconstruction error (or maximum performance metric, where applicable) and minimum cost, when compared to random sensor arrays.

Note that in the case of completely inaccessible locations (equivalent to infinite cost), the unmodified QR decomposition can be performed on the subset of allowed sensor locations.

\section{Balanced sensor and actuator selection for control systems}
\label{sec:control}

There is a large body of work regarding optimal sensor and actuator selection for control. Some, including~\cite{dhingra2014admm, munz2014sensor}, pose sensor and actuator selection as a convex optimization problem. Many take advantage of submodular or supermodular performance metrics to design greedy algorithms, for example~\cite{summers2014optimal, cortesi2014submodularity, summers2015submodularity, tzoumas2015minimal}. Nestorovi{\'c} and Trajkov~\cite{nestorovic2013optimal} design a sensor and actuator selection method for flexible structures such as cantilever beams based on balanced truncation and the $H_2$ and $H_\infty$ norms, which is similar to the framework considered here. We follow the method of Manohar et al.~\cite{manohar2018optimal} and perform the cost-constrained column-pivoted QR decomposition on truncated balanced modes for sensor and actuator selection. An outline of this method follows,~\cite{manohar2018optimal,Brunton2019book} contain full details.

Consider the control system
\begin{subequations}
\begin{align}
\dot{{\bf x}} &= {\bf A}{\bf x} + {\bf B}{\bf u}\\
{\bf y} &= {\bf C}{\bf x}.
\end{align}
\end{subequations}
We seek a suitable basis to place sensors (select rows of the identity matrix to form ${\bf C}$) and actuators (select columns of the identity matrix to form ${\bf B}$). We employ balanced truncation~\cite{moore1981principal, willcox2002balanced, rowley2005model}, which uses the controllability and observability Gramians to construct a basis in which the system is simultaneously maximally controllable and observable.

The controllability and observability Gramians ${\bf W}_c$ and ${\bf W}_o$ are given by
\begin{align}
{\bf W}_c = \int_0^\infty e^{{\bf A}t}{\bf B}{\bf B}^* e^{{\bf A}^*t} dt\\[4pt]
{\bf W}_o = \int_0^\infty e^{{\bf A^*}t}{\bf C}^*{\bf C} e^{{\bf A}t} dt.
\end{align}
Their eigenvectors define the directions in which the system is most controllable and observable, ranked by the magnitude of the corresponding eigenvalues. The balancing transformation makes the controllability and observability Gramians equal to each other and diagonal:
\begin{align}
\tilde{{\bf W}}_c = \tilde{{\bf W}}_o = {\bf\Sigma},
\end{align}
where
\begin{align}
\tilde{{\bf W}}_c &= {\bf\Phi}^*{\bf W}_c{\bf\Phi},\\
\tilde{{\bf W}}_o &= {\bf\Psi}^*{\bf W}_o{\bf\Psi},
\end{align}
${\bf\Phi}^* = {\bf\Psi}^{-1}$, and ${\bf\Sigma}$ is diagonal. The transformation matrix can be calculated by solving an eigenvalue problem:
\begin{align}
{\bf W}_c{\bf W}_o{\bf\Psi} = {\bf\Psi}{\bf\Sigma}^2.
\end{align}
With the eigenvalues ranked in decreasing order, the first $r$ modes of ${\bf\Psi}$ and ${\bf\Phi}$ define the $r$ most jointly observable and controllable directions for the system.

We use these modes to perform a rank-reduced change of basis, such that the new control system is written as:
\begin{subequations}
\begin{align}
\dot{{\bf a}}_r &= {\bf \Phi}^*_r{\bf A}{\bf \Psi}_r{\bf a}_r + {\bf \Phi}^*_r{\bf B}{\bf u}\\
{\bf y} &= {\bf C}{\bf \Psi}_r{\bf a}_r,
\end{align}
\end{subequations}
where ${\bf\Psi}_r$ and ${\bf\Phi}_r$ are the first $r$ columns of ${\bf\Psi}$ and ${\bf\Phi}$.

To select sensors and actuators with a cost function, we initialize ${\bf B}$ and ${\bf C}$ to identity matrices, calculate the Gramians and the truncated balancing modes, and perform the cost-constrained column-pivoted QR decomposition of Algorithm \ref{alg:ccqr} on ${\bf \Psi}^*_r$ for sensor selection and on ${\bf \Phi}_r^*$ for actuator selection.

As a performance metric, take the $H_2$ norm,
\begin{align}
\text{tr} \left({\bf C}{\bf W}_c{\bf C}^*\right) = \text{tr} \left({\bf B}^*{\bf W}_o{\bf B}\right),
\end{align}
which is a system's output energy, or the root mean square of its impulse response. A small output energy indicates that the sensor and actuator arrays are able to quickly stabilize the system after an input is applied, and therefore we wish to minimize the $H_2$ norm. However, we calculate a proxy:
\begin{align}
\log\det {\bf C}{\bf W}_c{\bf C}^* = \log\det {\bf B}^*{\bf W}_o{\bf B}.
\label{eq:H2proxy}
\end{align}
As shown in Summers~\cite{summers2015submodularity}, maximizing Eq. \ref{eq:H2proxy} is equivalent to minimizing the $H_2$ norm. Furthermore, Manohar~\cite{manohar2018optimal} demonstrates that $\arg\max_{\bf C} \log\det{\bf C}{\bf W}_c{\bf C}^*$ simplifies to $\arg\max_{\bf C}|\det{\bf C}{\bf\Psi}_r|$, hence why the greedy determinant maximization of the QR algorithm is desirable.

Our modification of adding a cost function during sensor and actuator selection removes any performance guarantees, but as will be shown in the next section, it is possible to effectively select sensors and actuators while accounting for a nonuniform cost on location.

\subsection{Example}
\label{subsec:ControlEx}

\begin{figure}[b]
    \centering
    \includegraphics[width = \columnwidth]{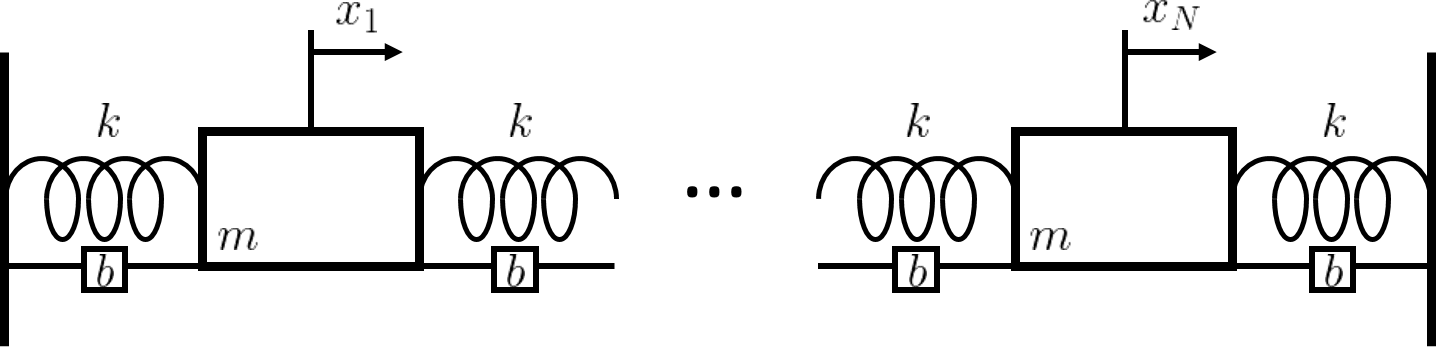}
    \vspace{-.25in}
    \caption{A one-dimensional system of $N$ identical masses with mass $m$, connected by springs with spring constant $k$, and damped by a factor of $b$.}
    \label{fig:MassSpring}
\end{figure}

We demonstrate balanced sensor and actuator selection on a system of $N$ identical masses with mass $m$, connected by identical damped springs, with spring constant $k$ and damping $b$~\cite{dhingra2014admm}, as illustrated in Fig. \ref{fig:MassSpring}. This is a simple benchmark system, and by choosing $N=16$ masses, it is possible to compare the principled sensor and actuator selections to a brute force search over every possible sensor/actuator configuration.

The equation of motion for the $i^{\text{th}}$ mass is given by
\begin{align*}
    m\ddot{x}_i = k(x_{i-1} + x_{i+1} - 2x_i) + b(\dot{x}_{i-1} + \dot{x}_{i+1} - 2\dot{x}_i) + f_i,
\end{align*}
where $x_i$ is the displacement of the $i^{\text{th}}$ mass from equilibrium and $f_i$ is an external forcing. The full state is ${{\bf x} = \left( \begin{array}{cccccc} x_1& \cdots& x_N& \dot{x}_i& \cdots& \dot{x}_N \end{array}\right)^T}$, and we let ${m = k = b = 1}$. Thus, the system with full actuation and sensing is given by
\begin{subequations}
\begin{align}
\label{eq:MassSpringDynamics}
    \dot{\bf x} &= \left(\begin{array}{cc} 
    {\bf 0}_N & \mathbb{I}_N\\[4pt]
    {\bf T} & {\bf T} 
    \end{array}\right) {\bf x} +
    \left(\begin{array}{c}
    {\bf 0}_N\\[4pt]
    \mathbb{I}_N
    \end{array}\right) {\bf f}\\[5pt]
    {\bf y} &= \mathbb{I}_{2N} {\bf x},
\end{align}
\end{subequations}
where ${\bf 0}_N$ is an $N\times N$ matrix of zeros and ${\bf T}$ is an $N\times N$ Toeplitz matrix with $-2$ along the diagonal and $1$ along the first super- and subdiagonals.

\begin{figure}[t]
\centering
\includegraphics[width=\columnwidth]{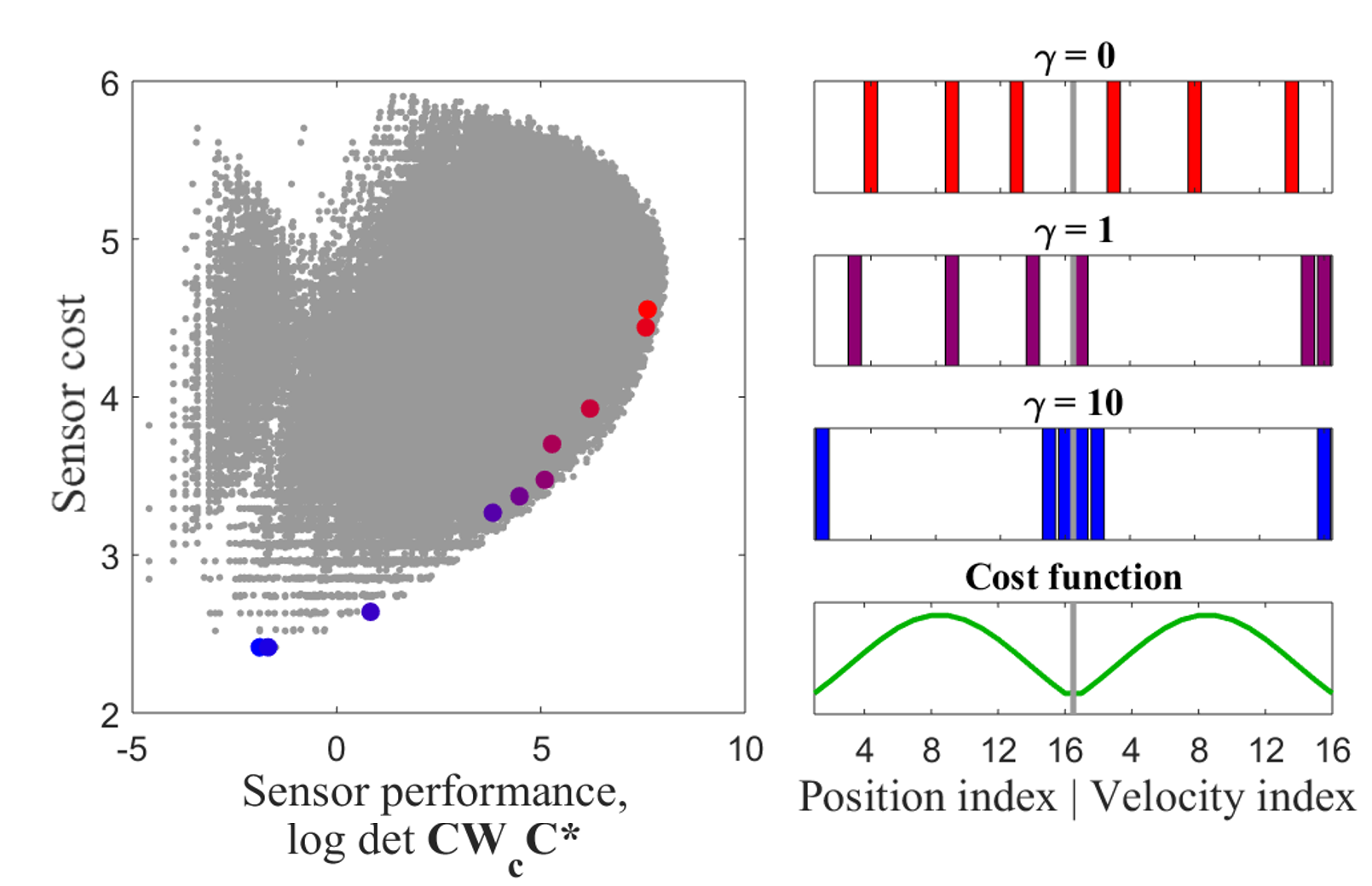}
\vspace{-.275in}
\caption{Sensor performance results and locations for the spring-mass system, placing 6 out of 32 possible sensors. The left plot shows cost versus performance metric log det ${\bf C}{\bf W}_c{\bf C}^*$. The large colored data points show the modified-QR-selected sensors as the cost function weighting is increased, with red representing $\gamma = 0$ and blue representing $\gamma = 10$. The gray points give the cost and performance metric for every possible configuration of sensor selection. On the right are the sensor locations at three values of $\gamma$, shown in the same colors as the corresponding points on the left. The Gaussian cost function is shown on the bottom right. The gray bar indicates the separation between position and velocity variables in ${\bf x}$.}
\label{fig:MassSpringResults_Sensors}
\vspace{-.1in}
\end{figure}

\begin{figure}[t]
\centering
\includegraphics[width=\columnwidth]{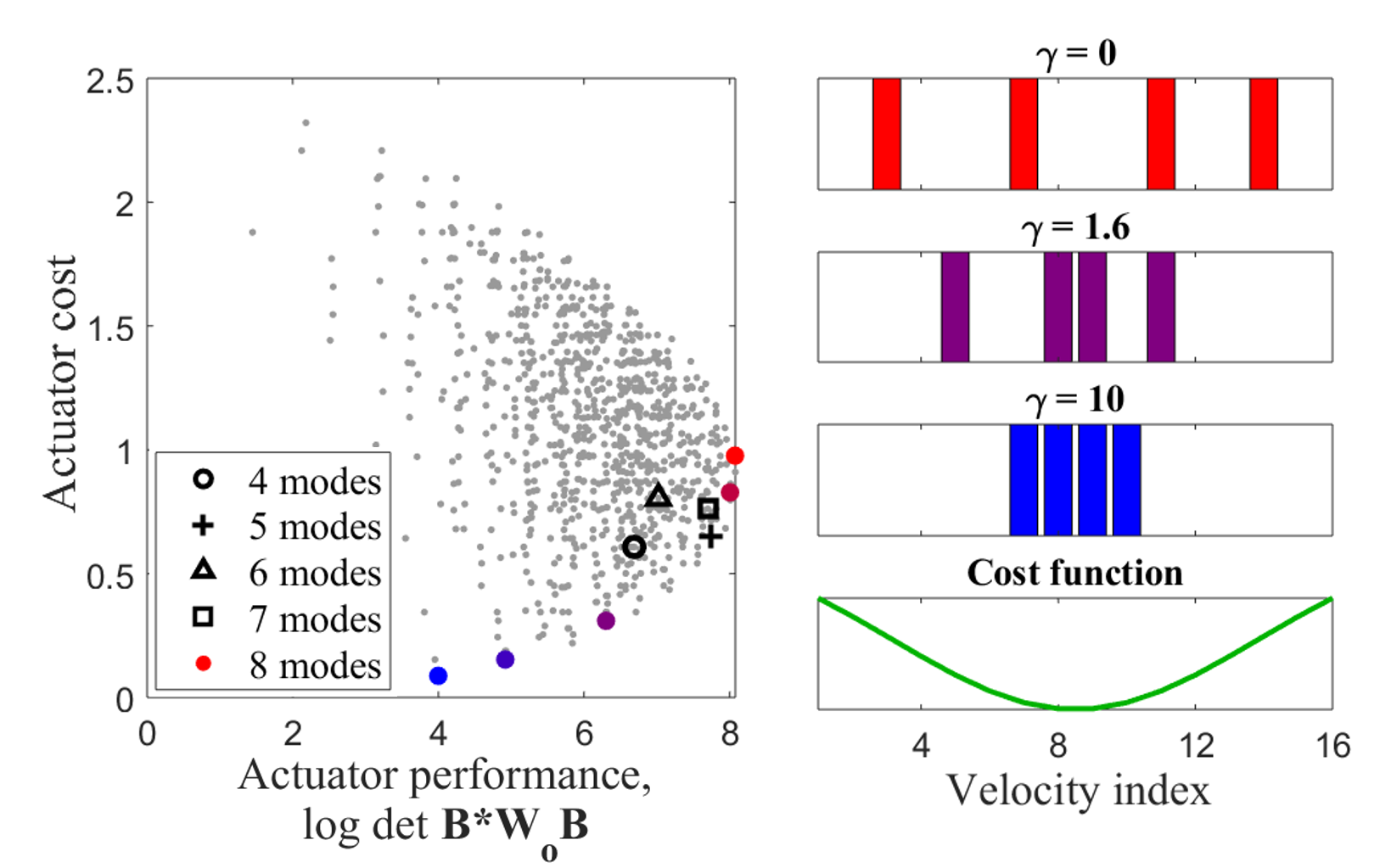}
\vspace{-.275in}
\caption{Actuator performance results and locations, placing 4 out of 16 possible actuators for the spring-mass system. Plotted on the left are cost versus performance metric log det ${\bf B}^*{\bf W}_o{\bf B}$. As in Fig. \ref{fig:MassSpringResults_Sensors}, the large colored data points show the results with principled actuator selections, shaded from red to blue as the cost function weighting is increased from zero to ten, while the gray points give the cost and performance metric for every possible configuration of actuator selection. The black symbols show the cost and performance using four to seven balanced modes with $\gamma = 0$. The right plots show three sets of actuator locations with colors corresponding to the plot on the left. The inverted Gaussian cost function is shown on the bottom right.}
\label{fig:MassSpringResults_Actuators}
\vspace{-.1in}
\end{figure}

We take $N = 16$ masses and place six sensors and four actuators. We choose a Gaussian cost function for the sensors, so is most expensive to place both position and velocity sensors in the center of the array. We select an inverted Gaussian as the actuator cost function, making it more expensive to place actuators at the edges of the array. The cost functions are shown in the lower right-hand panels of Fig. \ref{fig:MassSpringResults_Sensors} and \ref{fig:MassSpringResults_Actuators}.

Because the actuation takes the form of an externally applied force ${\bf f}$, which causes a change in the velocity of the masses, the actuators are restricted to acting on the velocity coordinates, which can be seen in Eq. \ref{eq:MassSpringDynamics}. Therefore, after calculating the truncated balanced and adjoint modes, we perform the cost-modified column-pivoted QR decomposition only on the final $N$ rows of the adjoint mode matrix for actuator selection. Equivalently, we could set the first $N$ rows to zero before performing the QR decomposition, or assign those coordinates an infinite cost. Because this cuts the domain of possible actuators in half, we find that we require twice the number of modes to achieve good performance, i.e. to place four actuators, truncate at eight balanced adjoint modes before performing the QR decomposition.

Results in terms of the $H_2$-norm metric substitute are shown in Fig. \ref{fig:MassSpringResults_Sensors} and \ref{fig:MassSpringResults_Actuators}, which plot cost versus $\text{log det }{\bf C}{\bf W}_c{\bf C}^*$ and $\text{log det }{\bf B}^*{\bf W}_o{\bf B}$ for the sensor and actuator arrays, respectively. The large colored data points are the results from the arrays chosen by QR pivoting, with red corresponding to $\gamma = 0$, and blue indicating a large value of $\gamma$. The gray points show the results from every other combination of choosing 6 out of 32 sensors or 4 out of 16 actuators. The principled sensors approximately follow the curve of simultaneously maximal $H_2$ proxy and minimal cost. The best-performing configuration (cost function weighting set to zero) outperforms 99.8\% of all other arrays, and with a high cost function weighting, the algorithm identifies the configuration with the lowest possible cost. Similarly, when the cost function weighting is set to zero, the selected actuators outperform all but three of the possible actuator selection permutations, and with a high value of $\gamma$, the algorithm again selects the lowest-cost permutation.

The plots on the right-hand sides of Fig. \ref{fig:MassSpringResults_Sensors} and \ref{fig:MassSpringResults_Actuators} show the sensor and actuator configurations for three different values of $\gamma$, in colors corresponding to those in the left panels. As the cost function weighting is increased, the sensors and actuators are pushed out of regions of high cost. The cost functions are also shown in green for reference. Finally, the left plot of Fig. \ref{fig:MassSpringResults_Actuators} also shows the performance of QR-selected actuators with $\gamma = 0$ and a number of balanced modes increasing from four to seven (shown in black, with symbols indicating the number of modes). Using eight modes yields the best results in terms of $\text{log det }{\bf B}^*{\bf W}_o{\bf B}$.

\begin{figure}[t]
\vspace{-.175in}
\includegraphics[width=\columnwidth]{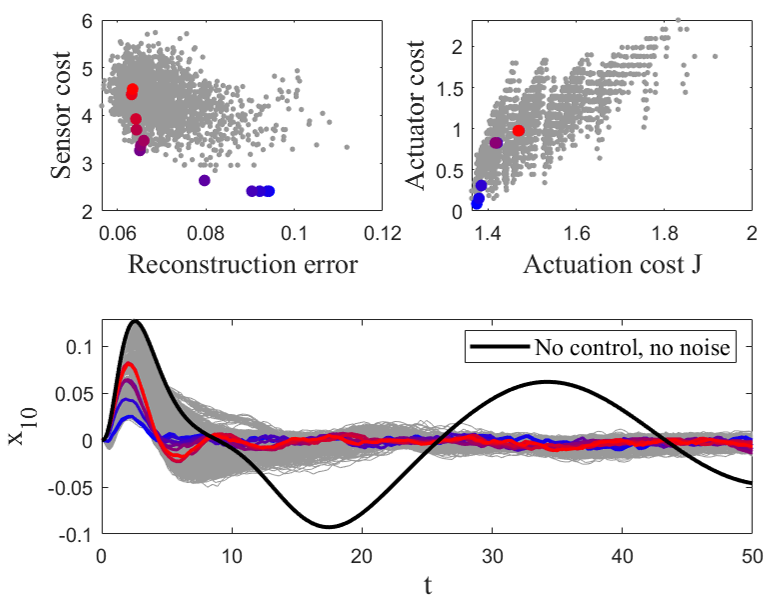}
\vspace{-.275in}
\caption{Sensor and actuator results using LQG control and state estimation for an example time series of the spring-mass system. The state is initialized at ${\bf x}(0) = 0$, except for $x_8(0) = 1$. The upper left plot shows sensor cost versus full-state reconstruction error. Red indicates $\gamma = 0$, blue corresponds to $\gamma = 10$, and gray shows results from 1820 random permutations of selecting 6 out of 32 possible sensors. The plot on the upper right similarly shows actuator cost versus the LQR actuation cost $J$. The gray dots show all 1820 possible configurations of placing 4 out of 16 actuators. Finally, we plot the trajectory $x_{10}(t)$ for the various sensor/actuator selections, with the black reference trajectory using no control and having no additive system or measurement noise.}
\label{fig:LQG_x8}
\vspace{-.125in}
\end{figure}

We also test sensor and actuator performance by evolving the system in time and applying linear-quadratic-Gaussian (LQG) control. To evaluate sensor performance, we add in a small amount of system disturbance and noise, both with covariance 0.005, and construct a Kalman filter to obtain a state estimation ${\bf \hat{x}}(t)$. We integrate the system from time $t=0$ to $t=50$ and stack the true states and their estimations into the columns of data matrices as shown:
\begin{subequations}
\begin{align}
{\bf X} &= \left[\begin{array}{cccc}
{\bf x}(0) & {\bf x}(t_1) & \cdots & {\bf x}(50)
\end{array}\right]\\[4pt]
{\bf \hat{X}} &= \left[\begin{array}{cccc}
{\bf \hat{x}}(0) & {\bf \hat{x}}(t_1) & \cdots & {\bf \hat{x}}(50)
\end{array}\right].
\end{align}
\end{subequations}
The reconstruction error is then given by Eq. \ref{eq:error}. Better sensor arrays will provide more informative measurements, and will therefore have lower reconstruction errors.

At every time step, the state estimate is passed to a linear-quadratic regulator (LQR), designed to drive the entire system towards ${\bf x} = {\bf 0}$. The performance of the actuator array is evaluated by calculating the LQR cost over the full timespan, 
\begin{equation}
J = \int_0^{50} \left({\bf x}^T {\bf Q}{\bf x} +  {\bf u}^T{\bf R}{\bf u}\right) dt,
\end{equation}
where ${\bf Q}$ and ${\bf R}$ are identity matrices. Good actuator configurations will be able to bring the system to equilibrium  quickly, and will have lower values of $J$.

\begin{figure}[t]
\centering
\vspace{-.175in}
\includegraphics[width=\columnwidth]{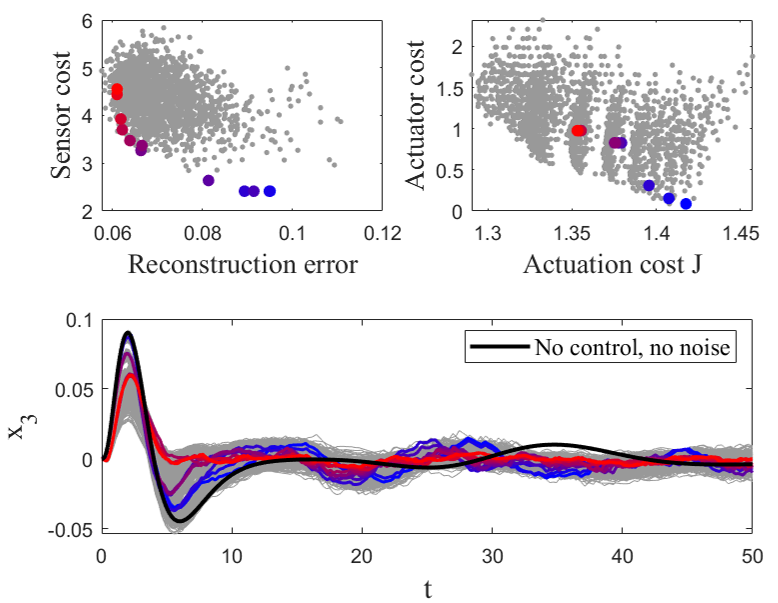}
\vspace{-.275in}
\caption{Sensor and actuator results for the spring-mass system using LQG control, with the inital condition ${\bf x}(0)=0$, except for $x_1(0) = 1$. As in Fig. \ref{fig:LQG_x8}, the upper left plot shows sensor cost versus reconstruction error, the upper right plot gives actuator cost versus $J$, and the lower plot gives an example time series, in this case $x_3(t)$. Red points and lines indicate a cost function weighting of zero, blue indicates a high value of $\gamma$, and gray gives non-principled sensor and actuator selections. The black line in the lower plot shows $x_3(t)$ with no control or noise.}
\label{fig:LQG_x1}
\vspace{-.125in}
\end{figure}

We test two different initial conditions, one with all positions and velocities set to zero except for the eighth mass, which has position $x_8(0)=1$, and the other where $x_1(0)=1$ with all other variables set to zero. Results are given in Fig. \ref{fig:LQG_x8} and \ref{fig:LQG_x1}. The figures have the same format, with the upper left subplot showing sensor cost versus reconstruction error, the upper right plot showing actuator price as a function of actuation cost $J$, and the lower subplot plotting the trajectory of one of the masses over time ($x_{10}(t)$ in Fig. \ref{fig:LQG_x8} and $x_3(t)$ in Fig. \ref{fig:LQG_x1}). The colored dots and lines show the results with the principled sensor and actuator arrays, with the usual color scheme of red indicating $\gamma=0$, shading through to blue for large values of $\gamma$. The gray points and lines are results from all 1820 permutations for placing 4 out of 16 actuators, each using a random array of sensors. The results for the upper plots are averaged over 25 realizations of noise, while the trajectories were all generated with the same noise vectors.

For both sets of initial conditions, the sensor configurations with $\gamma = 0$ provide effective, but not optimal, reconstructions compared to the randomly-selected arrays. As $\gamma$ is increased, the cost decreases and the reconstruction error goes up by about $4\%$ total. This can be seen qualitatively on the trajectory plots, where after $t\approx 20$, the red lines cluster more tightly about zero, while the blue lines have more variation.

As for the actuators, in the case of $x_1(0) = 1$, the QR-selected actuators are surprisingly far from optimal, though the performance decreases with the actuator cost, as expected. However, when $x_8(0)=1$, the actuator array chosen with the largest value of $\gamma$ leads not only to the lowest actuator cost but also to the optimum value of $J$. This seems counterintuitive, but because the center mass was given a large initial displacement, actuators placed in the center of the array are best suited to driving the entire system to zero in this case.

This example demonstrates that the column-pivoted QR decomposition on balanced modes for control systems can be modified to include a cost function on sensors and actuators. The method will select sensors and actuators to approximately simultaneously minimize cost and maximize performance, though the actual performance of the controller designed with the selected actuators may depend on the initial conditions.

\section{Sensor selection in a DMD basis}
\label{sec:DMD}
Dynamic mode decomposition (DMD) was introduced in the fluid dynamics community as a method of dimensionality reduction~\cite{schmid2010dynamic}. Conceptually, DMD is comparable to a spatial SVD combined with a temporal Fourier decomposition, in that it extracts spatially coherent modes that exhibit oscillatory or growth/decay behavior in time. This provides a highly effective reduced order model that can be used for short-time prediction, and can be adapted for applications such as control~\cite{proctor2016dynamic} and regime identification~\cite{kramer2017sparse}. Furthermore, DMD has strong connections to the Koopman operator~\cite{rowley2009spectral,Mezic2013arfm}: DMD approximates the Koopman operator advancing the state forward in time. Systems are often sparse in a DMD basis~\cite{jovanovic2014sparsity}, and DMD has been paired with compressed sensing in both space~\cite{brunton2015compressed,gueniat2015dynamic} and time~\cite{tu2014spectral}. Sparse sensing for multiresolution DMD was formulated by Manohar et al.~\cite{manohar2019optimized}. Here we provide a sparse sensing method for DMD modes allowing for a heterogeneous cost on sensor location.

We calculate the dynamic mode decomposition as follows~\cite{tu2014dynamic}: We begin by collecting a set of $m$ full-state snapshots, spaced equidistant in time (an alternative formulation does not require snapshots to be spaced equally in time~\cite{askham2018variable}),
\begin{align}
{\bf X} = \left[\begin{array}{cccc}
{\bf x}(t_1) & {\bf x}(t_2) & \cdots & {\bf x}(t_m) 
\end{array}\right].
\end{align}
We then organize the snapshots into two matrices offset from each other by one time step:
\begin{subequations}
\begin{align}
{\bf X}_1 &= \left[\begin{array}{cccc}
{\bf x}(t_1) & {\bf x}(t_2) & \cdots & {\bf x}(t_{m-1}) 
\end{array}\right]\\
{\bf X}_2 &= \left[\begin{array}{cccc}
{\bf x}(t_2) & {\bf x}(t_3) & \cdots & {\bf x}(t_{m}) 
\end{array}\right].
\end{align}
\end{subequations}
The goal of DMD is to find the best-fit matrix ${\bf A}$ such that 
\begin{align}
{\bf X}_2 \approx {\bf A}{\bf X}_1,
\label{eq:DMDFormula}
\end{align} i.e. the matrix ${\bf A}$ that advances the system forward by one time step. We take the $r$-mode truncated SVD of ${\bf X}_1$:
\begin{align}
{\bf X}_1 = {\bf U}_r {\bf\Sigma}_r{\bf V}_r^*.
\end{align}
The full least squares solution of Eq.~\ref{eq:DMDFormula} is given by 
\begin{align}
{\bf A} = {\bf X}_2 {\bf X}_1^\dagger = {\bf X}_2 {\bf V}_r{\bf\Sigma}_r^{-1}{\bf U}_r^*,
\end{align}
but since a primary goal of DMD is rank truncation, only the leading $r$ eigenvalues and eigenvectors of ${\bf A}$ are of interest here. Therefore, we calculate the projection ${\bf\tilde{A}}$ of ${\bf A}$ onto the leading $r$ spatial SVD modes,
\begin{subequations}
\begin{align}
{\bf\tilde{A}} &= {\bf U}_r^* {\bf A} {\bf U}_r\\[4pt]
&= {\bf U}_r^* {\bf X}_2 {\bf V}_r {\bf \Sigma}_r^{-1}.
\end{align}
\end{subequations}
We take the eigenvalue decomposition:
\begin{align}
{\bf\tilde{A}}{\bf W} &= {\bf W}{\bf \Lambda}.
\end{align}
The eigenvalues along the diagonal of ${\bf \Lambda}$ are identical to the leading $r$ eigenvalues of the full matrix ${\bf A}$, and the leading eigenvectors are given by
\begin{align}
{\bf\Psi} &= {\bf X}_2{\bf V}_r {\bf\Sigma}_r^{-1}{\bf W}.
\end{align}
The $k^{\text{th}}$ state ${\bf x}_k = {\bf x}(t_k)$ can then be expressed in terms of these eigenmodes and eigenvalues:
\begin{align}
{\bf x}_k &= {\bf \Psi}{\bf\Lambda}^{k-1}{\bf b},
\end{align}
or the continuous-time equivalent
\begin{align}
\label{eq:DMD_continuous}
{\bf x}(t) &= {\bf \Psi}e^{{\bf\Omega}t}{\bf b},
\end{align}
where ${\bf b}$ is the vector of amplitudes ${\bf b} = {\bf\Psi}^\dagger {\bf x}_1$, and ${\bf\Omega}$ is a diagonal matrix whose diagonal entries are the continuous-time eigenvalues $\omega_j = \log \lambda_j/dt$.

Thus, the eigenvectors ${\bf\Psi}$ are coherent spatial modes that form a reduced basis for the system of interest. We use this basis for sensor selection by performing the cost-modified column-pivoted QR decomposition on ${\bf\Psi}^*$.

\subsection{Example}
\label{subsec:DMDEx}

\begin{figure}[t]
\centering
\includegraphics[width=0.8\columnwidth]{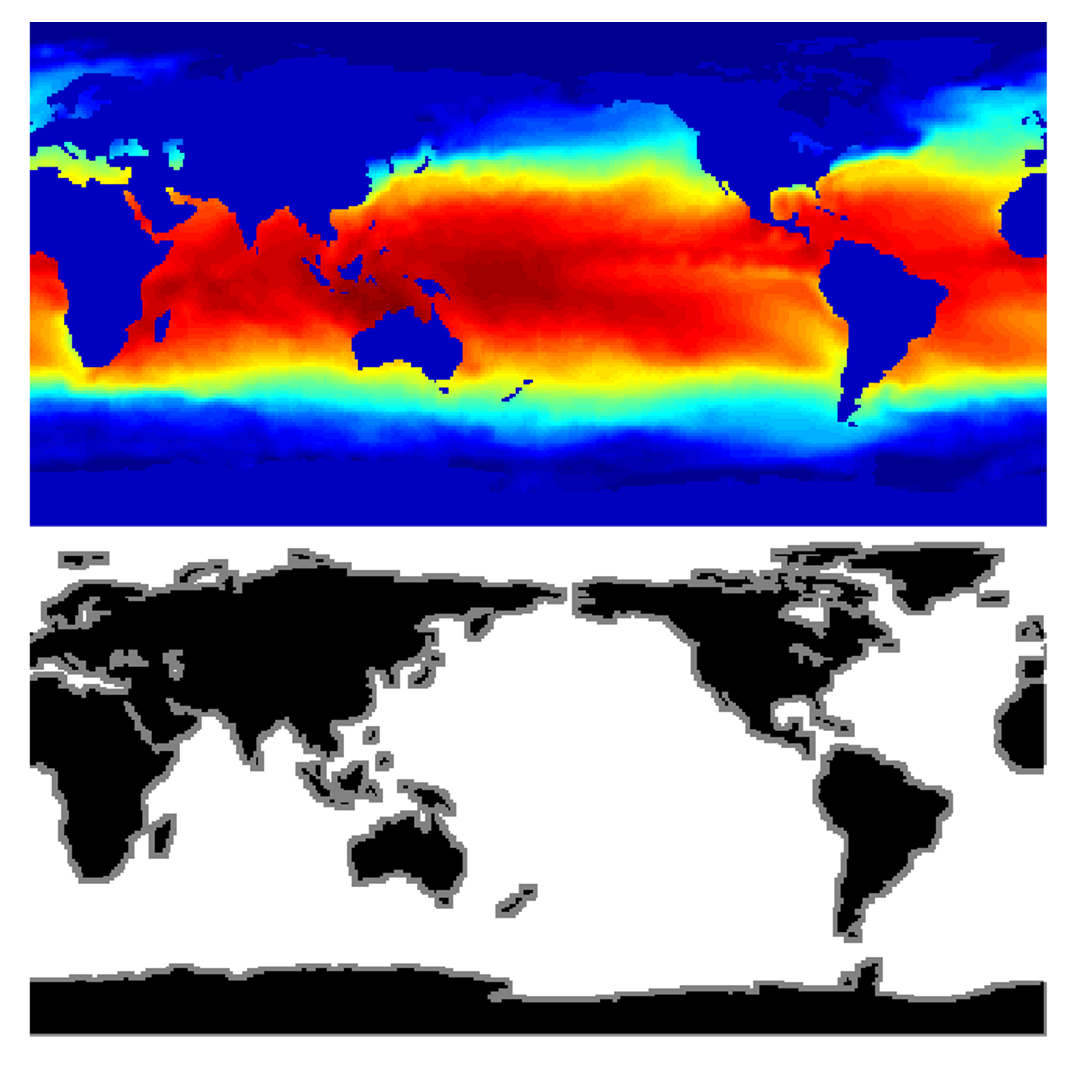}
\vspace{-.15in}
\caption{An example sea surface temperature snapshot (upper) and the cost function used (lower), where gray indicates regions of uniformly low cost and white denotes uniformly high cost.}
\label{fig:SSTSnapshotCost}
\end{figure}

As an example system, consider the NOAA weekly sea surface temperature measurements from 1990\textendash2016~\cite{noaa2018, banzon2016, reynolds2007}. This comprises 1455 temporal snapshots of sea surface temperature on a $360\times180$ spatial grid. The system is approximately periodic due to seasonal variations, so we expect DMD to give a reasonably good description of its behavior. To simulate a realistic sensor selection scenario, we use a cost function that uniformly penalizes placing sensors more than two grid points from land. An example sea surface temperature snapshot and the cost function are shown in Fig. \ref{fig:SSTSnapshotCost}.

We provide two performance metrics for the sensor arrays, interpolative error and extrapolative error:
\begin{subequations}
\begin{align}
E_{int} &= \frac{||{\bf X}_{tr}-{\bf\hat{X}}_{tr}||_{F}}{||{\bf X}_{tr}||_{F}}\\[4pt]
E_{ext} &= \frac{||{\bf X}_{te}-{\bf\hat{X}}_{te}||_{F}}{||{\bf X}_{te}||_{F}},
\end{align}
\end{subequations}
where the columns of ${\bf X}_{tr}$ are the training set, the first 80\% of snapshots, which are used to train the DMD modes. The columns of ${\bf\hat{X}}_{tr}$ are their reconstructions. ${\bf X}_{te}$ and ${\bf\hat{X}}_{te}$ are the test set, the remaining 291 snapshots and reconstructions sequentially following the training set.

\begin{figure}[t]
\centering
\vspace{-.1in}
\includegraphics[width = \columnwidth]{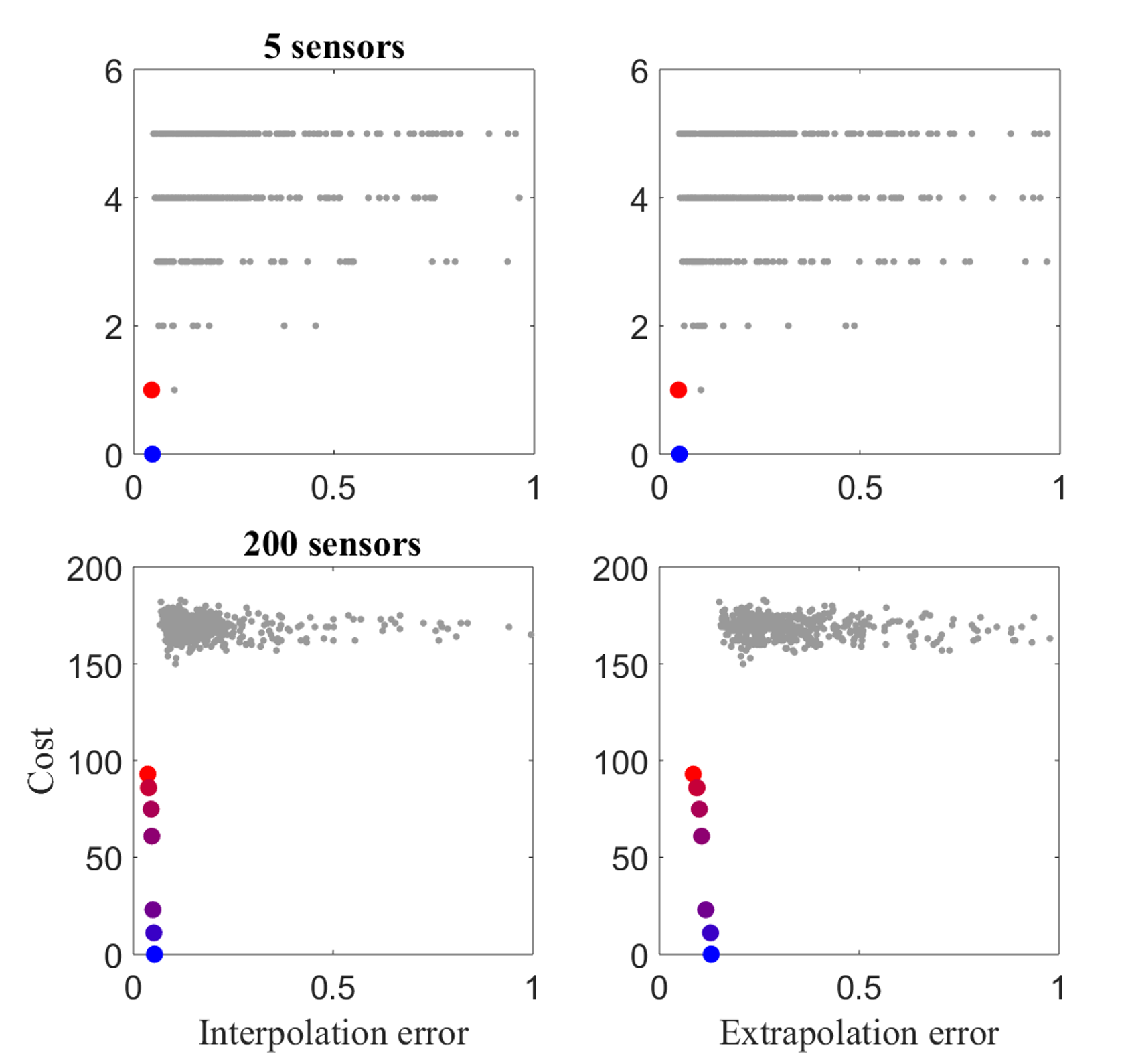}
\vspace{-.3in}
\caption{Cost versus reconstruction error for the sea surface temperature system using DMD modes as the cost function weighting is increased from zero (red) to one (blue). Results are shown for 5 sensors on the top row and for 200 sensors on the bottom row. The small gray points are results from 500 random sensor arrays.}
\label{fig:DMD_CostVsError}
\end{figure}

\begin{figure}
\centering
\vspace{-.1in}
\includegraphics[width = \columnwidth]{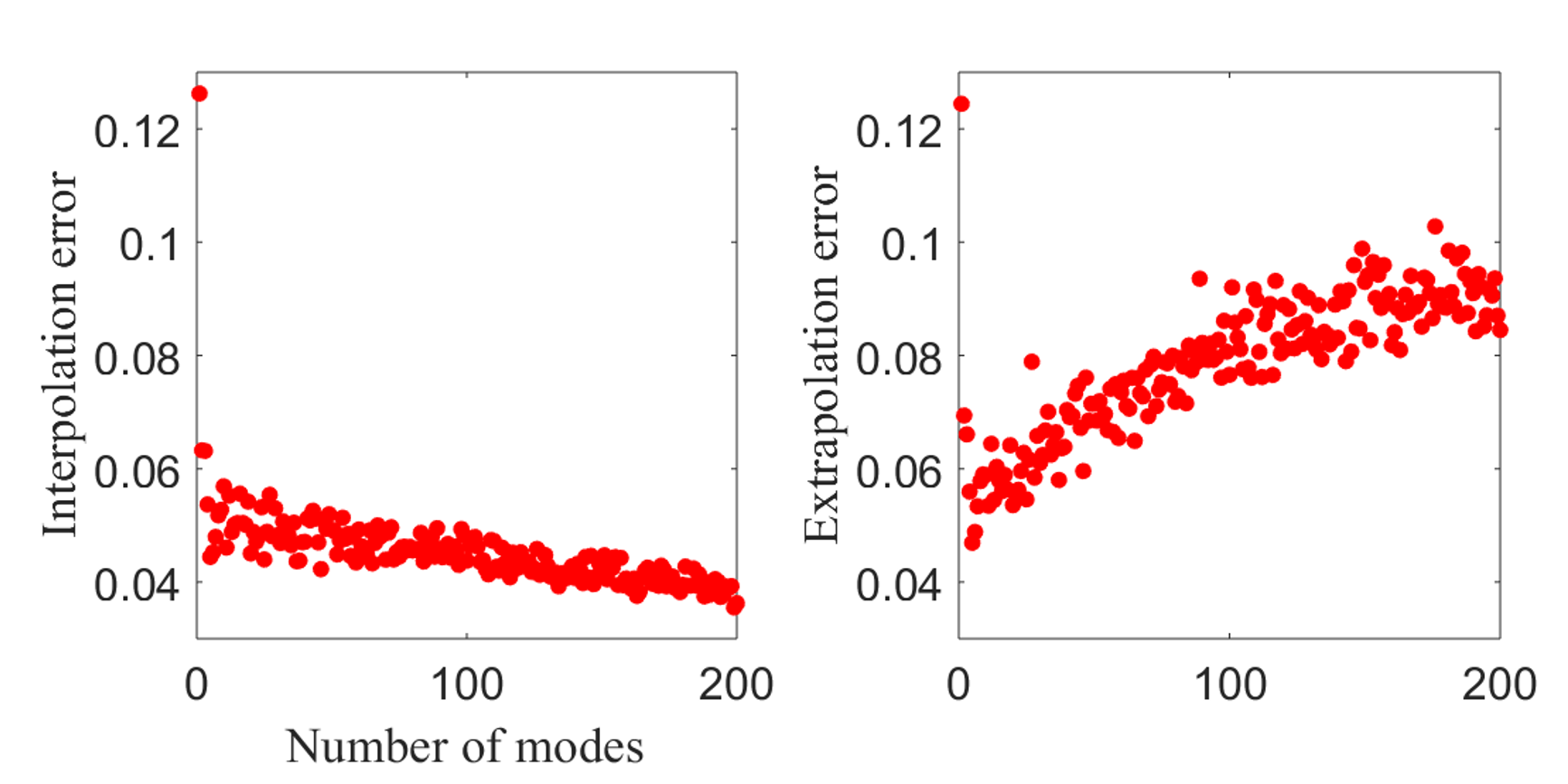}
\vspace{-.3in}
\caption{Error as a function of the number of DMD modes (equal to the number of sensors) for the sea surface temperature data set. Interpolation error is shown on the left and extrapolation error is on the right. No cost function is used.}
\label{fig:DMD_ErrorVsModes}
\end{figure}

The sea surface temperature data set has a slow singular value decay with a Gavish-Donoho cutoff~\cite{gavish2014optimal} around 400 modes. Thus we test our sensor selection method with both a small number of sensors (5) and a large number of sensors (200), though both cases represent significant downsampling from the original snapshots with dimension $\sim10^{5}$. Our results are shown in Fig. \ref{fig:DMD_CostVsError}. Red data points have a cost function weighting of zero, while blue points have a large weighting, $\gamma = 1$. As the cost function is weighted more heavily, the sensor array cost decreases to zero, leading to a small increase in reconstruction error, though remarkably this increase is by less than 1\%. Results from 500 random sensor arrays are also shown for comparison in gray. The QR-based algorithm places most of the sensors close to land even without incorporating the cost function, so not only do the random sensors generally lead to higher reconstruction errors, they also have higher costs. As expected, interpolation errors are lower than extrapolation errors. Unexpectedly, there is little improvement in increasing from 5 to 200 sensors, and in fact extrapolation errors are higher with more sensors. This is most likely related to overfitting and noise amplification, see~\cite{peherstorfer2018stabilizing} for more details. We plot interpolation and extrapolation error with no cost function for $r$ between 1 and 200 in Fig. \ref{fig:DMD_ErrorVsModes}. It is apparent that while interpolation error decreases slightly as the number of modes is increased, the extrapolation error increases significantly, as suggested by Fig. \ref{fig:DMD_CostVsError}.

\begin{figure}[t]
\centering
\vspace{-.1in}
\includegraphics[width = \columnwidth]{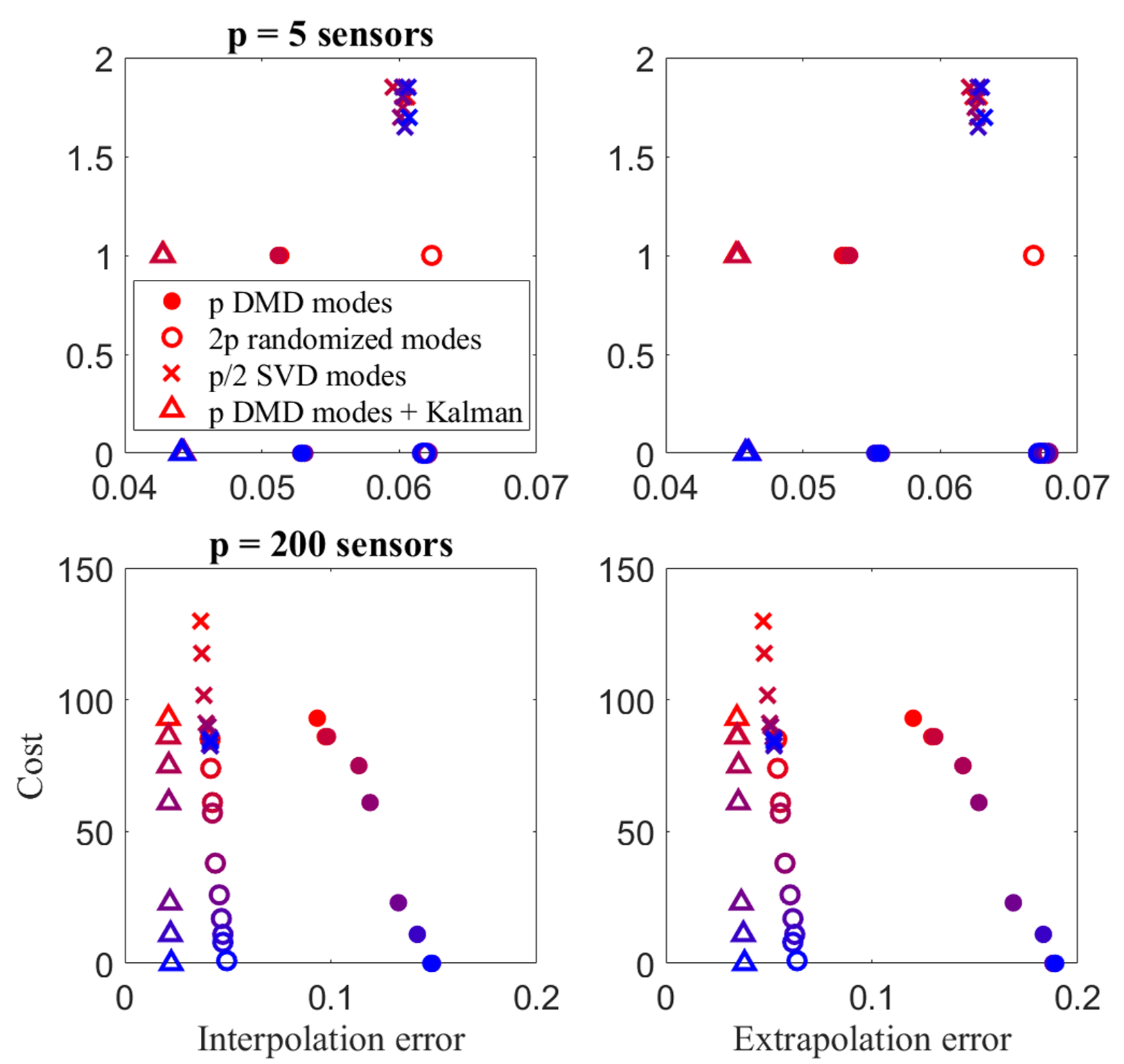}
\vspace{-.25in}
\caption{Cost versus reconstruction error for the sea surface temperature system, using DMD, randomized, and SVD modes with $p=5$ and 200 sensors. For the DMD modes, we test both the usual linear algebra reconstruction method of Eq. \ref{eq:recona} and a Kalman estimator (filled circles).}
\label{fig:DMD_CostVsError_Comp}
\end{figure}

We also provide a comparison to other basis choices in Fig.~\ref{fig:DMD_CostVsError_Comp}. We test randomized rank reduction as described in~\cite{halko2011finding}, and take twice the number of modes as sensors for greatly improved performance. These results are plotted as open circles. Results with an SVD basis are shown as crosses, and we follow the guidelines of~\cite{peherstorfer2018stabilizing}, performing the modified QR decomposition on the first $p/2$ SVD modes and placing the remaining $p/2$ sensors randomly to avoid noise amplification. Both of these basis choices lead to lower reconstruction errors at 200 sensors than the DMD basis, shown as closed circles, though the random sensors used with the SVD basis lead to higher costs. However, the benefit of DMD is the information it provides about the system's dynamics, which allows us to construct a Kalman estimator for the DMD amplitudes. From Eq. \ref{eq:DMD_continuous}, the time-varying amplitudes are given by
\begin{align}
{\bf b}(t) = e^{{\bf\Omega}t}{\bf b}(0).
\end{align}
We can then evolve the system forward in time according to
\begin{align}
\frac{d}{dt}{\bf {b}} = {\bf \Omega}{\bf b}.
\end{align}
We include added sensor noise with $2\%$ variance (this noise is applied to all of the sensor arrays for Fig. \ref{fig:DMD_CostVsError_Comp}). The Kalman estimator, results plotted as triangles, significantly outperforms the other three methods, suggesting that sensors and reconstruction methods that account for the dynamics of a system are preferable when possible.

\section{Sensor selection in an analytic basis}
\label{sec:analytic}
In many physical systems, a tailored data-driven basis is not necessary because the system is well-described by a universal basis, or because its equations of motion are known and can be solved analytically to yield a modal decomposition.  It is then possible to write the state according to Eq.~\ref{Eq:Representation}, where 
${\bf a}$ are Fourier coefficients and ${\bf \Psi}$ is a DFT matrix with $r=n$. Then the sparse sensing problem is similar to compressed sensing~\cite{candes2006robust}, in which a universal basis with random measurements is used to reconstruct a signal. Here, however, the signal is not required to be sparse in ${\bf \Psi}$, since we use least squares, rather than $L_1$-regularized reconstruction. Other examples of Eq.~\ref{Eq:Representation} include polynomial interpolation and Gaussian quadrature.

Alternatively, the system can be time-varying:
\begin{align}
{\bf x}(t) = {\bf\Psi}{\bf a}(t).
\end{align}
This is mathematically similar to the DMD formulation of Eq. \ref{eq:DMD_continuous}, the main differences being that the spatial modes are here assumed to be found analytically, rather than through data-driven methods, and we make no assumptions about the form of the time dependence.

As in previous sections, once ${\bf\Psi}$ is calculated, Algorithm \ref{alg:ccqr} can be applied to ${\bf\Psi}^*$ for cost-constrained sensor selection.

\subsection{Example}
\label{subsec:VandermondeEx}

We test the performance of an analytic basis for sensor selection on a vibrating two-dimensional circular membrane, since the equations of motion for this system are known and can be solved for a modal decomposition. The membrane satisfies the following partial differential equation:
\begin{align}
\frac{\partial^2 u}{\partial t^2} = c^2\left( \frac{1}{r} \frac{\partial}{\partial r} \left( r \frac{\partial u}{\partial r} \right) + \frac{1}{r^2} \frac{\partial^2 u}{\partial \theta^2}\right),
\end{align}
with $\theta\in(-\pi, \pi]$, $r\in[0,a]$. We take the boundary and initial conditions
\begin{subequations}
\begin{align}
&u(a,\theta,t) = 0\\
&u(r,\theta,0) = \alpha(r,\theta)\\
&\frac{\partial u}{\partial t}(r,\theta,0) =0.
\end{align}
\end{subequations}

This system can be solved analytically through separation of variables. The solution is a Bessel series in $r$ and a Fourier series in $\theta$ and $t$:
\begin{multline}
u(r,\theta,t) = \sum_{n=1}^\infty A_{0n}J_0 \left(\sqrt{\lambda_{0n}}r\right) \cos \left(c\sqrt{\lambda_{0n}}t\right) +\\
\sum_{m=1}^\infty \sum_{n=1}^\infty \left( A_{mn}\cos m\theta + B_{mn} \sin m\theta\right) J_m \left(\sqrt{\lambda_{mn}}r\right) \\
\times\cos\left( c\sqrt{\lambda_{mn}} t\right),
\end{multline}
where $\lambda_{mn} = (z_{mn}/a)^2$, with $z_{mn}$ the $n^\text{th}$ zero of the $m^\text{th}$ Bessel function.

Thus, the basis modes of the system are $\cos (m\theta) J_m \left(\sqrt{\lambda_{mn}}r\right)$ and $\sin (m\theta) J_m \left(\sqrt{\lambda_{mn}}r\right)$, forming the columns of the basis. Practically, the basis matrix must be finite, so we truncate at $m=M$, $n=N$ modes:

\begin{multline}
{\bf \Psi}_{MN}  = \Big( \begin{array}{cccc}
J_0\left(\sqrt{\lambda_{01}}r\right) & \cdots & J_0\left(\sqrt{\lambda_{0N}}r\right) & \cos (\theta)\end{array}\\[4pt] 
\begin{array}{cccc}
\times J_1\left(\sqrt{\lambda_{11}}r\right) & \sin (\theta) J_1\left(\sqrt{\lambda_{11}}r\right) & \cdots &\cos (\theta) J_1\left(\sqrt{\lambda_{1N}}r\right)\end{array}\\[4pt]
\begin{array}{ccc}
\sin (\theta) J_1\left(\sqrt{\lambda_{1N}}r\right) & \cdots & \sin (M\theta) J_M\left(\sqrt{\lambda_{MN}}r\right)
\end{array}\Big).
\label{eq:BesselPsi}
\end{multline}

The entire system evolves as
\begin{align}
{\bf u}(t) &= {\bf\Psi}_{MN} \cos \left(c\sqrt{{\bf\Lambda}_{MN}} t\right) {\bf b}_{MN},
\end{align}
where ${\bf\Lambda}_{MN}$ is a diagonal matrix with entries $\lambda_{01},\dots,\lambda_{0N},\lambda_{11},\lambda_{11},\lambda_{12},\lambda_{12},\dots,\lambda_{1N},\lambda_{1N},\dots,\lambda_{MN},\lambda_{MN}$, and the cosine is taken elementwise on the diagonal entries. The vector ${\bf b}_{MN}$ gives the coefficients:
\begin{multline}
{\bf b}_{MN} = \big(\begin{array}{cccccccc}
A_{01} & \cdots & A_{0N} & A_{11} & B_{11} & \cdots & A_{1N} & B_{1N}
\end{array}\\[4pt]
\begin{array}{ccc}
\hspace{22pt}\cdots & A_{MN} & B_{MN}
\end{array}\big)^T.
\end{multline}

Our goal is to sparsely sample in space and then estimate the time-varying coefficients $\cos \left(c\sqrt{{\bf\Lambda}_{MN}} t\right) {\bf b}_{MN}$, subsequently obtaining a reconstruction of the full state ${\bf{\hat{u}}}(t)$. We select sensors by performing the cost-constrained column-pivoted QR decomposition on ${\bf\Psi}_{MN}^T$.

\begin{figure}[t]
\includegraphics[width=\columnwidth]{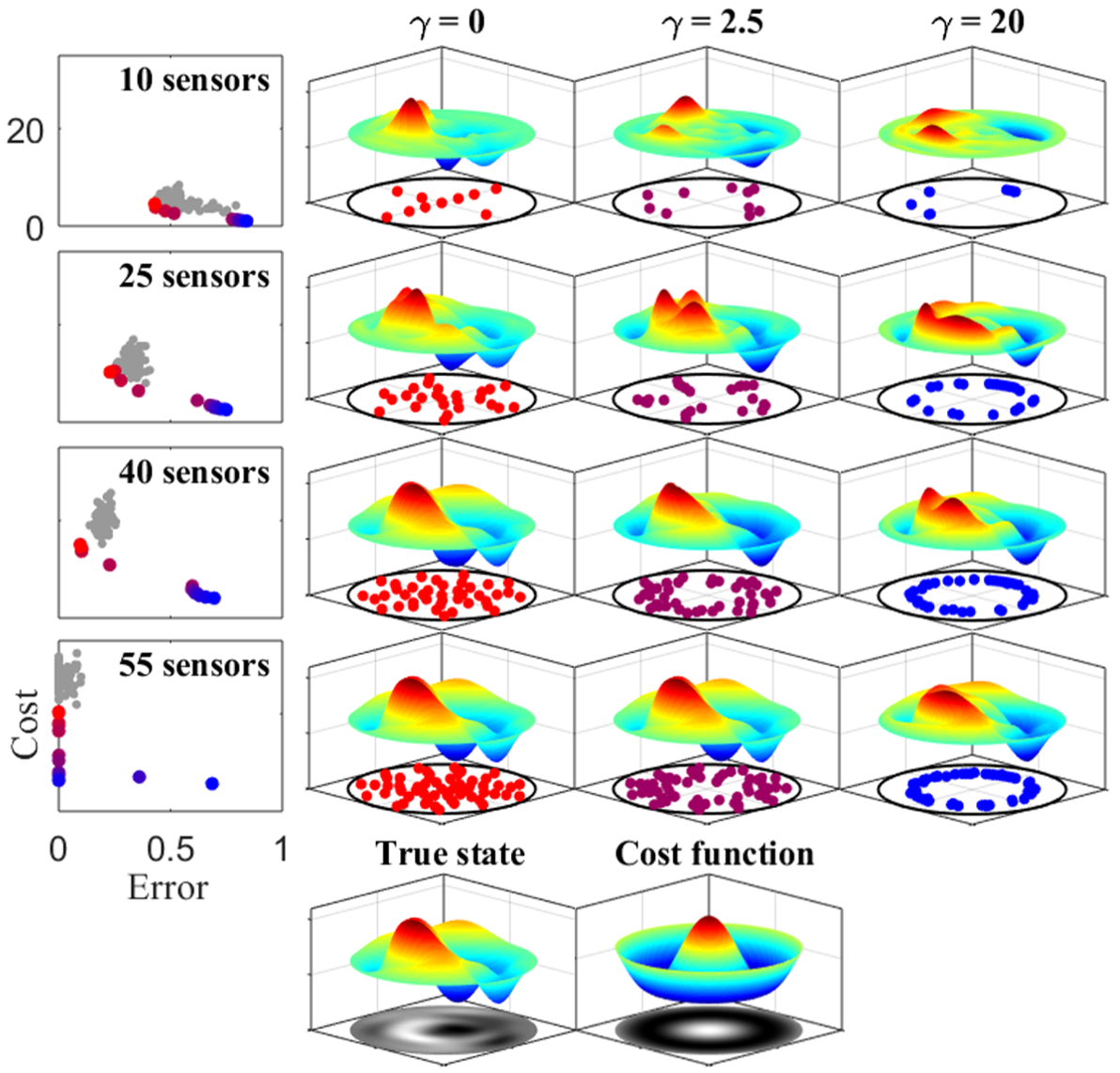}
\vspace{-.3in}
\caption{Results for sensor selection with a cost function for the vibrating drum system with 55 Bessel function basis modes. The number of sensors increases moving down the rows. The first column shows the cost versus reconstruction error as the cost function weighting is increased (red points indicate $\gamma=0$, blue points have high values of $\gamma$, gray points show 100 random sensor arrays). Note that all plots in the first column have the same $x$- and $y$-axis scales. The remaining three columns show reconstructions of an example snapshot with the indicated number of sensors and value of $\gamma$. The sensor locations are given as well, with color corresponding to that on the cost versus error plots. The true state and the cost function are shown on the bottom row, in both a surface and a pcolor plot.}
\label{fig:Bessel_CostError_Recons}
\end{figure}

For this experiment, we take $a = 10$, $c = 1$, and a grid of 101 points each in $r$ and $\theta$. We also truncate the modes at $M=6$ and $N=5$ for a total of 55 basis modes. The system is initialized with random coefficients, with lower modes weighted more heavily, i.e. 
\begin{align}
A_{mn}, B_{mn} \propto \frac{1.5}{n(m+1)}.
\label{eq:BesselCoeffs}
\end{align} Not only does this lead to more accurate reconstructions when undersampling, it is also often the case in physical systems that the lower modes are the most active. We choose a radially symmetric cost function with a maximum in the center and a minimum at $r=13/2$, $f(r) = 0.6 + 0.5\cos \left(\frac{2\pi r}{13}\right)$, pictured in the lower right panel of Fig. \ref{fig:Bessel_CostError_Recons}.

\begin{figure}[t]
\centering
\includegraphics[width = .95\columnwidth]{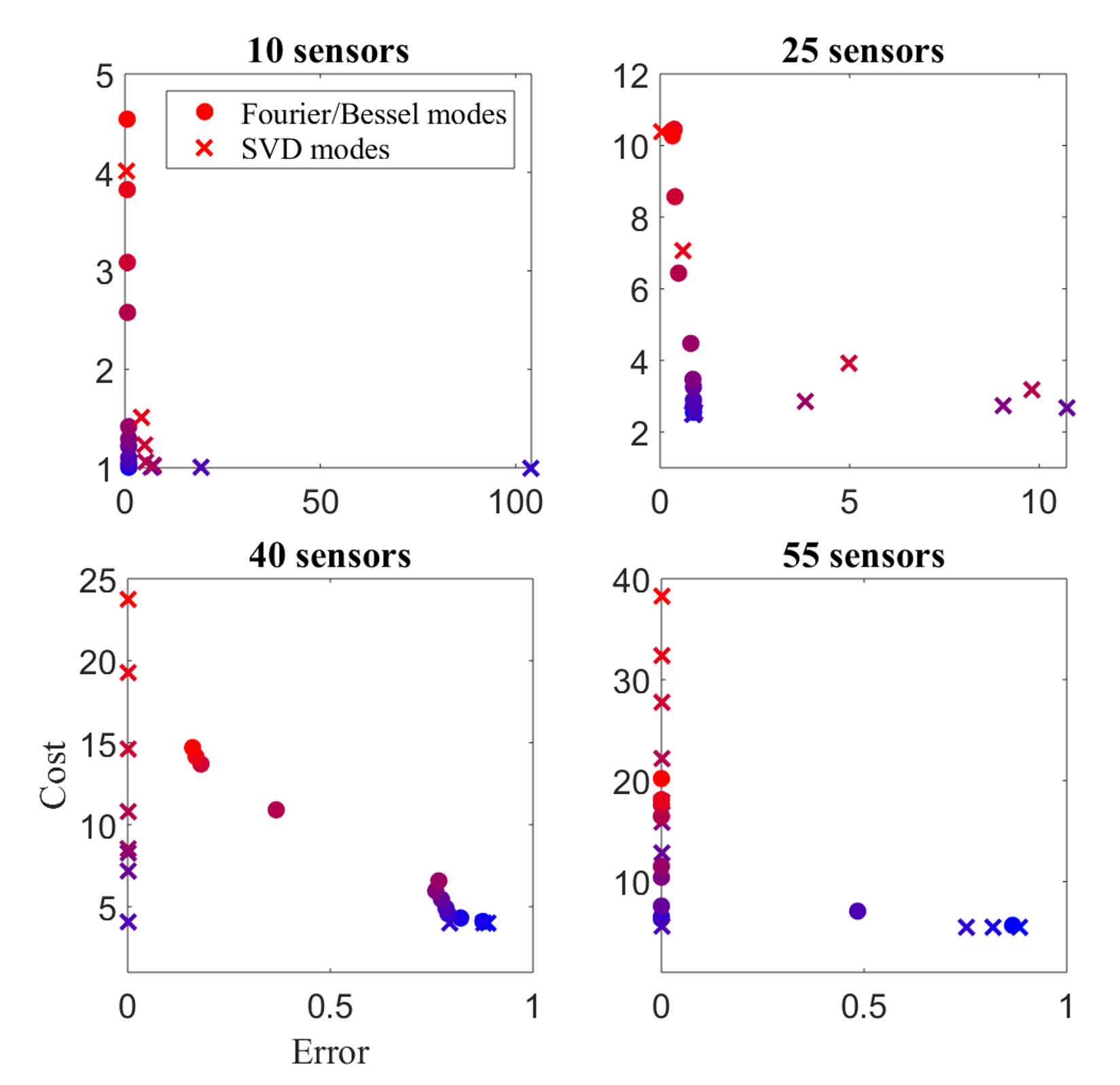}
\vspace{-.2in}
\caption{Cost versus reconstruction error for the vibrating drum system, using the analytic Fourier/Bessel modes of Eq. \ref{eq:BesselPsi} and SVD modes calculated from snapshots. As in Fig. \ref{fig:Bessel_CostError_Recons}, results are shown for different numbers of sensors and multiple cost function weightings, with the usual color scheme of red indicating $\gamma = 0$ shading through to blue representing a large value of $\gamma$.}
\label{fig:Bessel_CostVsError_Comp}
\end{figure}

We take $\gamma\in[0,20]$, and evaluate sensor performance by propagating the initial conditions forward in time steps of 0.1 up to $t=10$. The true states are stacked into the columns of a snapshot matrix ${\bf U}$, and the error is given by $||{\bf U}-{\bf \hat{U}}||_{F}/||{\bf U}||_{F}$, where ${\bf\hat{U}}$ is the reconstruction obtained with the chosen sensor array. We calculate this error with 50 different sets of random initial conditions and take the average, plotted in the first column of Fig. \ref{fig:Bessel_CostError_Recons}, which shows sensor cost plotted against this average error as $\gamma$ is increased (red data points correspond to $\gamma=0$, blue show $\gamma = 20$). The remaining panels in the figure show reconstructions of an example state with the sensor arrays chosen with various numbers of sensors and values of $\gamma$. The true state is shown on the bottom left of the figure. Even with just 10 sensors, when $\gamma=0$, the reconstruction is imperfect but retains some main features of the true state. With 55 sensors, the state can be reconstructed exactly. As the cost function weighting increases, the reconstruction quality degrades, even with 55 sensors. A large jump in reconstruction error occurs when there are no longer any sensors close to the center of the domain. The figure also shows results from 100 random sensor arrays, which suggest that the principled sensors follow the curve of minimum cost and reconstruction error quite well.

We also compare reconstruction results with sensors selected using an SVD basis, shown in Fig. \ref{fig:Bessel_CostVsError_Comp}, which plots cost versus reconstruction error for 10, 25, 40, and 55 sensors. Sensors selected with analytic Fourier/Bessel modes are shown as circles, with SVD results shown as crosses. In both cases, red indicates a cost function weighting of zero and blue corresponds to a large cost function weighting. To calculate the SVD basis, we initialize a random system with coefficients as in Eq. \ref{eq:BesselCoeffs} and evolve it forward for 1000 time steps to obtain snapshots. We then randomly place 700 snapshots into a training matrix and perform the SVD to obtain a basis. After using this basis and Algorithm \ref{alg:ccqr} to select sensors, we test them on the remaining 300 snapshots, and compare to reconstructions of the same snapshots using sensors from the analytic basis. We take the same number of SVD modes as sensors since there is no noise in the system, and take the average cost and error over 50 random training and test sets.

\begin{figure}[t]
\centering
\includegraphics[width = \columnwidth]{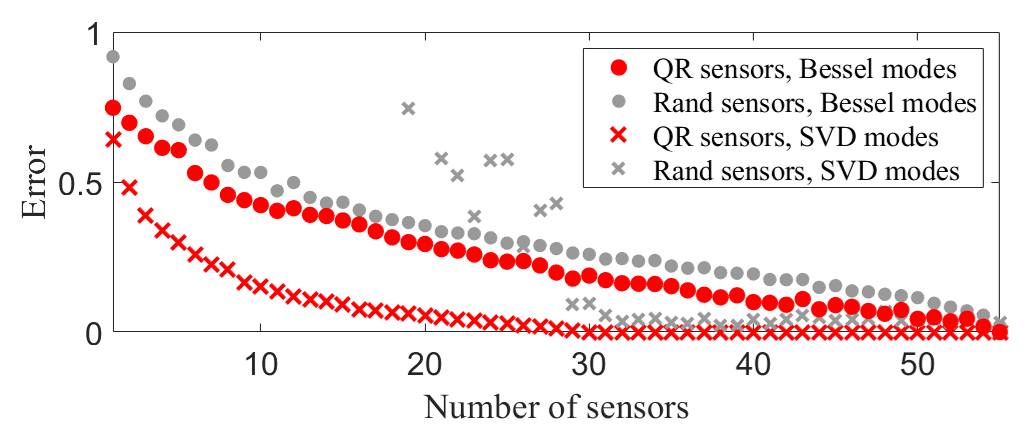}
\vspace{-.25in}
\caption{Error versus the number of sensors for the vibrating drum system, using both analytic Fourier/Bessel modes and data-driven SVD modes. Principled sensors selected by the QR algorithm are shown in red, compared to random sensors with both bases in gray. No cost function is used.}
\label{fig:Bessel_ErrorVsNum}
\end{figure}

As demonstrated in Fig. \ref{fig:Bessel_CostVsError_Comp}, while the SVD and analytic modes have comparable performance at a small number of sensors, the SVD basis is able to reconstruct the system perfectly to machine precision at 40 sensors. With the SVD basis, the active modes are mixed together, which is more beneficial for rank reduction\textemdash with sensors selected using analytic modes, even with 54 sensors there may be little information about the $55^{\text{th}}$ mode. This can be seen clearly in Fig. \ref{fig:Bessel_ErrorVsNum}, which shows error versus the number of sensors using both bases and no cost function. The SVD basis consistently outperforms the analytic basis, requiring just 30 sensors to reconstruct the full state perfectly, whereas the analytic modes need the full 55. However, the SVD requires data to calculate the basis, so in cases where high-fidelity full state data is not available, we have shown that analytically calculated modes are a viable alternative. Furthermore, the figure shows comparisons to random sensors using both bases for reconstruction. In the analytic basis, principled sensors have lower reconstruction errors than random sensors, and random sensors perform very poorly in an SVD basis, highlighting the importance of choosing the right method for the application at hand.

The results in this work indicate that if a dynamical system may be represented as the evolution of a superposition of known basis modes, these modes can be used for sparse sensor selection, even with a cost function on sensor location. Though a data-driven SVD basis leads to better sensor arrays, analytic modes often give comparable reconstructions. 
In a more complex signal, such as a sound signal being represented in a Fourier basis, where both high and low frequencies can be equally present and it is unknown which frequencies will be dominant, it is difficult to effectively downsample using this method. The problem is similar to that of compressed sensing, where in most real-world examples, the number of sensors required to reconstruct the signal with high probability is large and leads to a computationally prohibitive reconstruction problem. These examples demonstrate how data-driven methods are usually more efficient at extracting patterns from a system than historical analytic-function-based methods.

\section{Conclusions}
\label{sec:conclusion}

In this work, we explore cost-constrained QR pivoting for sensor selection, applying the algorithm to three new types of basis modes: balanced modes for sensor and actuator selection for control systems, DMD modes for dynamical systems, and universal or analytic modes. We find that all three types of modes make effective bases for sensor selection for their respective systems, and that the cost-constrained algorithm selects sensors that are nearly Pareto optimal in cost and performance. The new bases do not necessarily outperform a standard SVD basis, but they can still be beneficial. Sensors placed using a DMD basis can be used to construct a Kalman estimator, which does perform better than the SVD or randomized bases when measurement noise is present. And when no full-state data is available, analytic modes can be a good substitute for a data-driven basis.

These results suggest that sensor selection methods must account for real-world restrictions. A researcher or engineer can improve their full state estimations by selecting an appropriate basis, but they must also consider practicalities like a sensor cost landscape, hence the cost-constrained QR algorithm. There are also other practical considerations, such as noise and the cost of the sensors themselves. In this paper we added measurement noise only when constructing a Kalman filter, but physical sensors always have noise. Future work will examine sensor placement with multiple types of sensors, where expensive sensors have low noise levels and cheap sensors have higher measurement noise.

It will also be interesting to apply the cost-constrained QR algorithm to more complex systems, in particular fluid flows, such as the flow near a turbulent jet~\cite{arndt1997proper}, flow near an adjustable flap~\cite{taylor2004towards}, and control of the flow over an open cavity~\cite{cattafesta2003review, rowley2006dynamics}. Many interesting systems are multiscale in space or time, and we could consider sensor selection for such systems, similar to~\cite{manohar2019optimized}. It is clear that sensor selection is a complex topic, and despite continual advances, it is still open to further practical refinement for real-world applications.

\section*{Acknowledgment}
E. Clark and S. L. Brunton acknowledge support from the Boeing Company (664755). J. N. Kutz acknowledges support from the Air Force Office of Scientific Research (FA9550-19-1-0011).
S. L. Brunton acknowledges support from the Air Force Office of Scientific Research (FA9550-18-1-0200).

\ifCLASSOPTIONcaptionsoff
  \newpage
\fi







%

\begin{IEEEbiography}[{\includegraphics[width=1in]{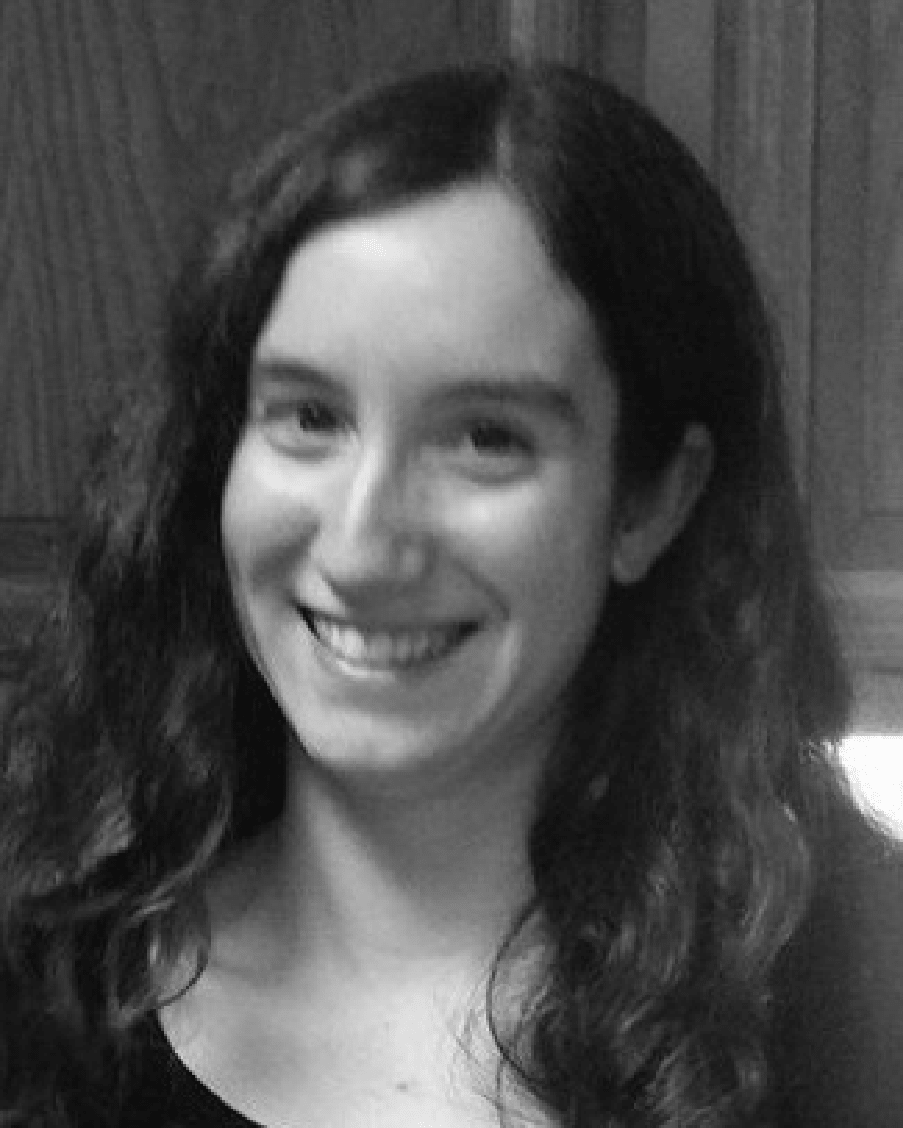}}]{Emily Clark}
received the B.S. degree in physics from Bates College, Lewiston, ME, in 2015. She is a Ph.D. candidate of physics at the University of Washington.
\end{IEEEbiography}

\begin{IEEEbiography}[{\includegraphics[width=1in]{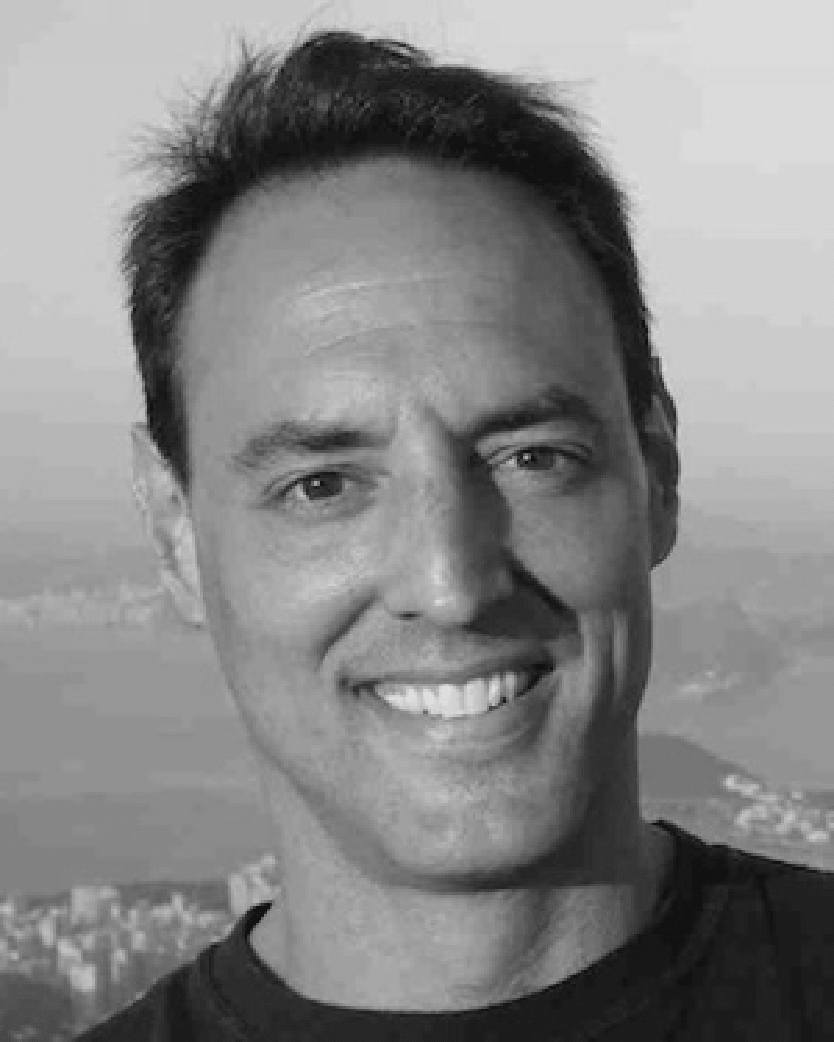}}]{J.~Nathan Kutz}
(Member, IEEE) received the B.S. degrees in physics and mathematics 
from the University of Washington, Seattle, WA, in 1990, 
and the Ph.D. degree in applied mathematics from 
Northwestern University, Evanston, IL, in 1994.  He is 
currently a Professor of applied mathematics, adjunct 
professor of physics and electrical engineering, and a 
senior data science fellow with the eScience institute 
at the University of Washington.
\end{IEEEbiography}

\begin{IEEEbiography}[{\includegraphics[width=1in]{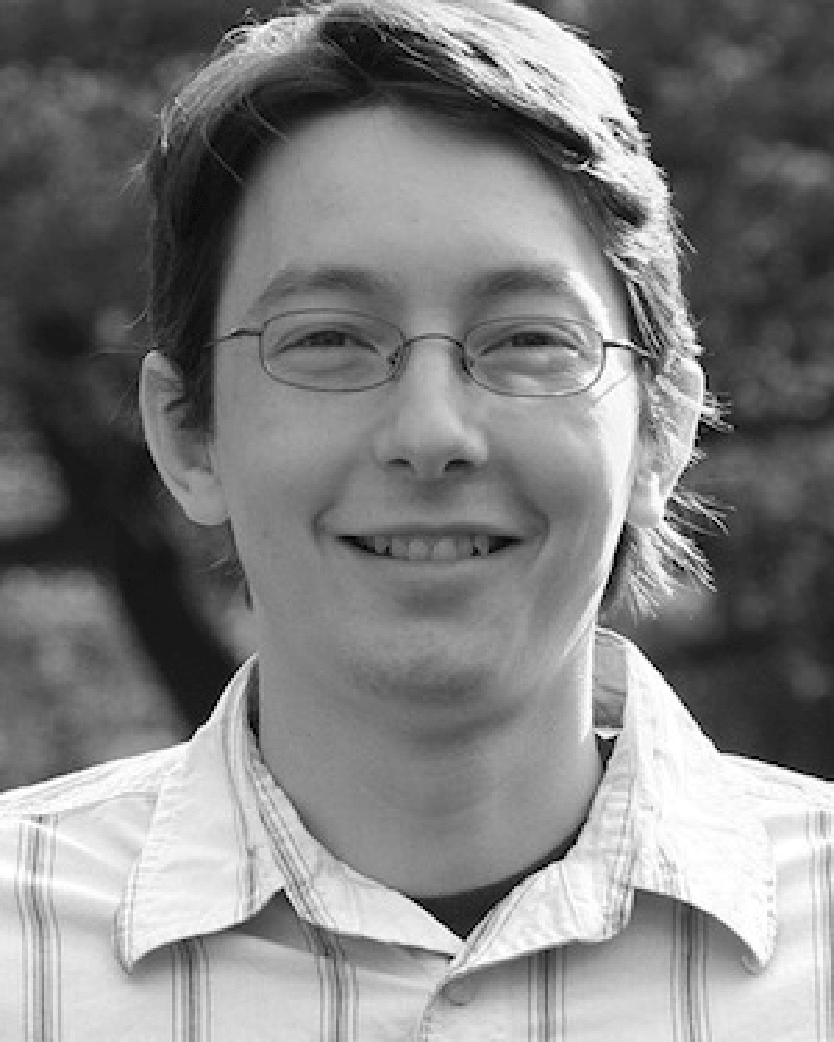}}]{Steven~L.~Brunton}
(Senior Member, IEEE) received the B.S. degree in mathematics with a minor 
in control and dynamical systems from the California 
Institute of Technology, Pasadena, CA, in 2006, and 
the Ph.D. degree in mechanical and aerospace engineering 
from Princeton University, Princeton, NJ, in 2012.  He 
is an Associate Professor of mechanical engineering and 
a data science fellow with the eScience institute at the 
University of Washington.
\end{IEEEbiography}

%







\end{document}